\documentclass{article}

\usepackage{lineno,hyperref}
\usepackage{amssymb}
\usepackage{amsmath}
\usepackage{amsfonts}

\usepackage{tikz-cd} 

\newcommand{\norm}[1]{\left\Vert#1\right\Vert}
\newcommand{\nnorm}[1]{{\left\vert\kern-0.25ex\left\vert\kern-0.25ex\left\vert #1 \right\vert\kern-0.25ex\right\vert\kern-0.25ex\right\vert}}

\newcommand{\field}[1]{\mathbb{#1}}
\newcommand{\R}{\field{R}}

\newcommand{\N}{\field{N}}
\newcommand{\M}{\field{M}}
\newcommand{\E}{\field{E}}
\newcommand{\F}{\field{F}}

\def\proof {\textbf{} }

\def\endproof{${}$\hfill$\triangle$}

\newtheorem{Theo}{Theorem}
\newtheorem{theo}{Theorem}[section]
\newtheorem{defi}[theo]{Definition}
\newtheorem{lem}[theo]{Lemma}

\newtheorem{prop}[theo]{Proposition}
\newtheorem{rem}[theo]{Remark}
\newtheorem{rem-a}{Remark}
\newtheorem{theo-a}[rem-a]{Theorem}
\newtheorem{nota}[theo]{Notations and conventions}
\newtheorem{com}[theo]{Comments}

\def\ap{\rightarrow}
\def\a{\alpha}
\def\b{\beta}
\def\g{\gamma}
\def\G{\Gamma}
\def\t{\tau}
\def\d{\delta}
\def\s{\sigma}
\def\so{\underline}
\def\dis{\displaystyle}
\def\p{\partial}

\begin{document}

\title{Regulated curves on a Banach manifold and singularities of endpoint map. I. Banach manifold structure}

\author{
T. Goli\'nski\footnote{Faculty of Mathematics, University of Białystok, email: tomaszg@math.uwb.edu.pl}
\and
F. Pelletier\footnote{UMR 5127 du CNRS \& Universit\'e de Savoie Mont Blanc, Laboratoire de math\'ematique, Campus scientifique, 73376 Le Bourget-du-Lac Cedex; 
 email: fernand.pelletier@univ-smb.fr}}

\maketitle

\begin{abstract}We consider regulated curves in a Banach bundle whose projection on the basis is continuous with regulated derivative. We build a Banach manifold structure on the set of such curves. This result was previously obtained for the case of strong Riemannian Banach manifold and absolutely continuous curves in \cite{Sch}. The essential argument used was the existence of a ``local addition'' on such a manifold. Our proof is true for any Banach manifold. 
In the second part of the paper the problems of controllability will be discussed.
\end{abstract}

%

\setcounter{tocdepth}{1}
\tableofcontents

\section{Introduction}
Classically, a sub-Riemannian structure on a finite dimensional manifold $M$ is a subbundle $E$ of $TM$ provided with a Riemannian metric $g$. An absolutely continuous path $\g: I := [a, b] \to M$ of class $H^1$ is called horizontal if the derivative $\dot{\g}$ is tangent to $E$ almost everywhere. Given any Riemannian metric $\hat{g}$ on $TM$ whose restriction to $E$ is $g$, it is well known that the set $\mathcal{H}^1(I, M)$ of paths $\g:[a,b]\to M$ of class $H^1$ has a structure of manifold and the set $\mathcal{H}^1(I, E)$ of horizontal curves is a closed Hilbert submanifold. Now given any point $x\in M$ the set $\mathcal{H}^1
_{x}=\{\g\in \mathcal{H}^1(I, M) \;: \g(a)=x\}$ is a finite codimensional submanifold and the map $\textrm{End}_{x}: \mathcal{H}^1_{x}\to M$ defined by $\textrm{End}_x(\g)=\g(b)$ is called the {\it endpoint map} which is smooth (cf. \cite{Mon}, \cite {AOP} and \cite{PiTa} for a more recent proof). If $\textrm{End}_{x}$ is surjective, then one says that the control system associated to this problem is called {\it controllable} and any critical point of $\textrm{End}_{x}$ is traditionally called an abnormal path. The existence of such curves is a big difference between Riemannian and sub-Riemannian case related to the problem of characterization of minimal geodesics and it was the origin of a lot of works about it (for nice references on this problem see \cite{Mon} or \cite{AgSa} for a more general context in control theory).

The concept of sub-Riemannian geometry in infinite dimensional context was introduced in \cite{GMV} for sub-bundle $\mathcal{D}$ of the tangent bundle $TM$ of a Riemannian ``convenient manifold'' $M$ with application to Fr\'echet Lie groups. Essentially motivated by {\it mathematical analysis of shapes}, the particular Banach context was recently studied by S. Arguill\`ere in his PhD thesis and more precisely formalized in \cite{Ar}. He considers a Banach bundle $\pi: E\ap M$ provided with an anchor $\rho:E\ap TM$ (called {\it anchored Banach bundle}) and a strong Riemannian metric $g$ on $E$. From the existence of a strong Riemannian metric it follows that the typical fiber of $E$ is a Hilbert space $\mathbb{H}$. Given any chart $(U,\phi)$ in $M$ such that $E_{| U}$ is trivial, a path $\g:I\to U$ is called {\it horizontal} or {\it admissible} if there exists a trivialization $\Phi: E_{| U}\to U\times \mathbb{H}$ and a path $u:I\to \mathbb{H}$ of class $L^2$ such that $\frac{d}{dt}(\phi\circ \g)=T\phi(\rho(\Phi^{-1}(\g,u))$ almost everywhere.
The set $\mathcal{H}^1(I, U)$ of a such curves has a Banach manifold structure since $\mathcal{H}^1(I,\phi(U))$ is open in the Sobolev space $H^1(I, \phi(U))\times L^2(I,
\mathbb{H})$ (cf. \cite{Ar} section 1). As in finite dimension, it follows that $\mathcal{H}^1_{x}=\{\g\in \mathcal{H}^1(I, M) \;: \g(a)=x\}$ is a closed submanifold and the 
map $\textrm{End}_{x}: \mathcal{H}^1_{x}\to M$ is again well defined and smooth. However, the range of $T_\g\textrm{End}_{x}(T_\g \mathcal{H}^1_{x})$ may not be closed but only dense in $T_{\textrm{End}_{x}(\g)}M$ which gives rise to another type of singular point $\g$ called {\it elusive path} in \cite{Ar}.

The main purpose of this work is to give a {\it generalization of such local results to global ones}. Thus in the context of an anchored bundle $(E,\pi, M,\rho)$, we need to show that a certain set of paths $\mathcal{P}(I, E)$ has a Banach manifold structure and the subset $\mathcal{AP}(I, M)$ of admissible paths has a structure of Banach submanifold of $\mathcal{P}(I,M)$. Moreover as in finite dimension, we want to show that for each $x\in M$ the endpoint map $\mathcal{AP}_{x}(I, M)\to M$ is smooth.

\bigskip

The first part of the paper deals only with the Banach manifold structure on $\mathcal{P}(I, E)$. {\it We now give a short outline of our results.}

\bigskip

For simplicity for now $I$ will denote the interval $[0,1]$. Let $M$ be a Banach manifold modelled on a Banach space $\mathbb{M}$. A path $c: I\ap M$ is called { \it regulated }
if for any $t\in I$ 
both one-sided limits exist. Such a path is called {\it strong regulated} or {\it $0$-regulated} if 
it is right-continuous everywhere and left-continuous at $1$ (cf. $\S$\ref{regulated}). A path $c: I\ap M$ is called { \it $1$-regulated } if $c$ is differentiable almost everywhere and its derivative is a $0$-regulated path \footnote{Such a path is absolutely continuous in the sense of \cite{Glo}.}.
For $k=0,1$ the set $\mathcal{R}^k(I,\mathbb{M})$ of $k$-regulated paths has a natural structure of Banach space (cf. Proposition \ref{Banachstr}).

Let $\pi:E\ap M$ be a Banach bundle with fibers modeled on Banach space $\E$. 
The set of $1$-regulated paths in $M$ defined on $I$ is denoted $\mathsf{P}(M)$ and the set of $0$-regulated paths $c:I\to E$ such that $\pi\circ c$ belongs to $\mathsf{P}(M)$ is denoted by $\mathcal{P}(E)$.

 
The main result of this part is:

\begin{Theo} \label{1}${}$
\begin{enumerate}
\item The set $\mathsf{P}(M)$ has a Banach manifold structure modeled on the Banach space $\mathcal{R}^1(I,\M)$. Moreover for any $\g\in \mathsf{P}(M)$ we have a Banach isomorphism from $T_\g\mathsf{P}(M)$ onto the set $\G(\g)=\{\mathsf{X}\in \mathcal R^1(I,TM)\;|\; p_M(\mathsf{X})=\g\}$.
\item The set ${\mathcal{P}}(E)$ has a Banach manifold structure modeled on the Banach space $\mathcal{R}^1(I,\M)\times\mathcal{R}^0(I ,\E)$. For any $c\in \mathcal{P}(E)$ we have a Banach isomorphism from $T_c {\mathcal{P}}(E)$ onto $\G(c)=\{\mathcal{X}\in \mathcal R^0(I,TE)\;|\; p_E(\mathcal{X})=c\}$. Moreover $\pi:\mathcal{P}(E)\ap \mathsf{P}(M)$ is a Banach bundle with typical fiber $\mathcal{R}^0(I ,\E)$.
\end{enumerate}
\end{Theo}
 
This Theorem can be obtained from \cite{Sch} for absolutely continuous paths when $M$ is provided with a strong Riemannian metric which implies the existence of a ``local addition'' on $M$ (cf. \cite{Mic}). Our proof goes through building an atlas on $\mathsf{P}(M)$ and $\mathcal{P}(E)$ which is a generalization of the proof for $H^1$ paths in sub-Riemannian finite dimension (cf. \cite{AOP}). In fact, it is also an adaptation of proofs of results of \cite{Pen1} or \cite{Kri} for maps of class $C^k$ from some compact manifolds with corners $\Omega\subset \mathbb{R}^m$ to a Banach manifold $M$. However, we can not deduce Assertion 1 directly from \cite{Kri} since the set $\mathcal{R}^0(I,\E)$ does not satisfy the axioms of section 2 in \cite{Kri}. Thus, to be complete, we give detailed proof of these results. Note that, according to the recent paper of Gl\"ockner \cite{Glo}, our proof can be transposed to the context of absolutely continuous paths as defined in this work.

This paper tries to be self contained. The next section \ref{regulated} contains essentially all the useful notations and results about $k$-regulated curves in Banach spaces. Some of these results were proved in \cite{Di,Pen2}. All needed results about $k$-regulated paths in the case of a Banach manifold are contained in $\S$\ref{Banachmanifold}. The fundamental results on the Banach structure on the set $\G(\g)$, $\G(c)$ and $\G(\g^*(E))$ (Proposition \ref{Banachstr}) can be found in $\S$\ref{GcGg}. The subsection \ref{topologyPE} describes the topology on $\mathcal{P}(E)$ and $\mathsf{P}(M)$. Theorem \ref{1} is proved in $\S$\ref{PEn}.
We have included Appendix A devoted
to the case of $C^1$-paths in Banach bundle $E$ such that $\pi\circ c$ is $C^2$-path. Such a result is basic for \cite{Pe}.
\section{$k$-regulated paths in a Banach space}
\label{regulated}

Consider a Banach space $\mathbb{E}$ equipped with a norm $\norm{\;}_\E$. Let $I=[t_0,t_1]$ be a compact interval. According to \cite{Di}, a path $c:I\ap \mathbb{E}$ is called {\it regulated} if for any $t\in[t_0,t_1[$ the limit $\dis\lim_{s\rightarrow t^+}c(s)$ exists and for any $t\in ]t_0,t_1]$ the limit $\dis\lim_{s\rightarrow t^-_1}c(s)$ exists. 
Note that the set of discontinuities of $c$ is at most countable and in particular it has {\it zero Lebesgue measure}. Now if $c:I\ap \mathbb{E}$ is regulated, the restriction of $c$ to any subinterval $J= [s_0,s_1]$ is a regulated path $c_J$ as well. 

\bigskip

\begin{defi}
A regulated path $c:I\ap \E$ is called {\it strong regulated} if it is right-continuous everywhere and left-continuous at $t_1$, i.e. for any point $t\in [t_0,t_1[$ we have $c(t)=\dis\lim_{s\rightarrow t^+}c(s)$ and $c(t_1)=\dis\lim_{s\rightarrow t_1^-}c(s)$.
\end{defi}

Note that for any regulated path $c:I\ap \E$, there exists a unique strong regulated path $\bar{c}$ 
that coincides with $c$ outside of discontinuity set and its discontinuity set is not larger than that of $c$. 

\bigskip 

If $c:[t_0,t_1]\ap \mathbb{E}$ and $c':[t'_0,t'_1]\ap \E$ are two strong regulated paths 
then {\it the concatenation $c{\star}c'$} is the strong regulated path defined by
\begin{align*}
c\star c'(t) &= c(t)\qquad&&\textrm{for } t\in [t_0,t_1[,\\
c\star c'(t_1) &= \dis\lim _{t\ap t_0'^+} c'(t),\\
c\star c'(t) &= c'(t-t_1)\qquad&&\textrm{for } t\in] t_1,t_1+t'_1-t'_0].
\end{align*}

The set $\mathcal{R}(I,\mathbb{E})$ of regulated paths in $\mathbb{E}$ provided with the norm $\norm{c}_{\infty}=\dis\sup_{t\in I}\norm{c(t)}_\E$ is a normed space. Let $U\subset\mathbb{E}$ be an open neighbourhood of the image $c(I)$ of a regulated path $c$. Then for any (local) diffeomorphism $\phi$ from $U$ to an open set $V$ it is clear that the path $\phi\circ c:I\ap \mathbb{E}$ is also regulated.

According to \cite[section VIII.7]{Di} for any $c\in \mathcal{R}(I,\mathbb{E})$ there exists a primitive of $c$ that is a continuous path $ \hat{c}:I\ap \mathbb{E}$ which is differentiable outside a countable set $\Sigma$ and such that $\dis\frac{d \hat{c}}{dt}=c(t)$ for all $t\in I\setminus \Sigma$. Moreover all such primitives are of type $\hat{c}+\operatorname{const}$. 
Note that if we impose the condition of strong regularity then $c$ is completely determined by any of its primitives. 


Note that while concatenation 
of two strong regulated paths 
(as defined previously) is a strong regulated path, the restriction to any subinterval $J\subset I$ of strong regulated path $c$ defined on $I$ {\it is not in general a strong regulated path on $J$} but only regulated. Therefore we consider: 

\begin{defi}\label{restrict}
Let $c:I\ap \E$ be a strong regulated path $c:I\ap \E$ and $J$ any subinterval of $I$. We define the strong regulated path $c_J$ as $\overline{c|_J}$. When there will be no risk of confusion we will call $c_J$ the restriction of $c$ to $J$.
\end{defi}
\begin{rem}\label{restrictstrong}\normalfont 
For strong regulated $c$ and $J=[s_0,s_1]\subset I$ the strong regulated path $c_J$ coincides with $c|_J$ on $[s_0,s_1[$ and its value at $s_1$ is $\dis\lim_{s\ap s_1^-}c(s)$.
\end{rem}

\begin{defi} For $k> 0$ a path $c:I\ap \mathbb{E}$ is called $k$-regulated 
if $c$ is of class $C^{k-1}$ and there exists a strong regulated path $c^{(k)}:I\ap \mathbb{E}$ such that the derivative of order $k-1$ of $c$ is a primitive of $c^{(k)}$. The strong regulated path $c^{(k)}$ is unique and will be called $k^{th}$ derivative of $c$.
\end{defi}
Note that if $c:I\ap \mathbb{E}$ is $1$-regulated, then we have 
\begin{equation}\label{c^1}
c(t)=c(t_0)+\dis\int_{t_0}^t c^{(1)}(s)ds. 
\end{equation}
Therefore for $k$-regulated curve $c$, for $1\leq l\leq k$, we always have:
\begin{equation}\label{c^l}
c^{(l-1)}(t)=c^{(l-1)}(t_0)+\dis\int_{t_0}^t c^{(l)}(s)ds.
\end{equation}

\begin{nota} \label{nota} ${}$\\ 
By $0$-regulated path we will mean a strong regulated path. For any $k \geq 0$ the set of $k$-regulated paths in $\E$ defined on an interval $I$ will be denoted $\mathcal{R}^k(I,\E)$. We provide $\mathcal{R}^k(I,\E)$ with the norm
$$\nnorm{c}_\infty^k=\dis\max_{0\leq l\leq k}\sup_{t\in I} \norm{c^{(l)}(t)}_\E.$$
\end{nota}

For an open set $U$ of $\E$ we can consider the subset
$$\mathcal{R}^k(I,U)=\{ c\in \mathcal{R}^k(I,\E)\;:\; c(I)\subset U\}.$$
Clearly $\mathcal{R}^k(I,U)$ is an open subset of $\mathcal{R}^k(I,\E)$.

\bigskip

We have the following classical results:

\begin{theo}\label{propreg} ${}$
\begin{enumerate}
\item For any $k$-regulated path $c:I\ap \E$ the set $\overline{c(I)}$ is compact.
\item The set $\mathcal{R}(I,\mathbb{E})$ provided with the norm $\norm{\;}_{\infty}$ is a Banach space.
\item The set $\mathcal{R}^k(I,\mathbb{E})$ for $k\geq0$ provided with the norm $\nnorm{\:}^k_\infty$ is a Banach space.
\item Let $A$ be a continuous linear map from $\E$ to a Banach space $\F$. If $c\in \mathcal{R}^k(I,\mathbb{E})$ then $A\circ c$ belongs to $ \mathcal{R}^k(I,\mathbb{F})$ and for any $1\leq l\leq k$ we have 
$$A\left(\dis\int_{t_0}^{t_1}c^{(l)}(s)ds\right)=\dis\int_{t_0}^{t_1}(A\circ c^{(l)})(s)ds.$$
\item Let $c:I\ap \E$ be a regulated curve and $\phi: J= [s_0,s_1]\ap \R$ be $1$-regulated curve such that $\phi(J)\subset I$. If either $c$ is continuous or $\phi$ is strictly monotone then $s\mapsto \phi'(s)(c\circ\phi)(s)$ is regulated and we have:
$$\dis\int_{s_0}^{s}\phi'(r)(c\circ\phi)(r)dr=\dis\int_{\phi(s_0)}^{\phi(s)}c(t)dt.$$
\end{enumerate}
\end{theo}

\noindent {\bf Proof.}\\
The fundamental observation is that the closure of the set of all step functions (resp. strong regulated step functions) on $I$ in the norm $\norm{\;}_{\infty}$ is $\mathcal{R}(I,\mathbb{E})$ (resp. $\mathcal{R}^0(I,\mathbb{E})$), see \cite[7.6.1]{Di} or \cite[Proposition 2.15]{Pen2}.
\begin{enumerate}
\item For $k=0$ the usual proof goes along these lines: we observe that since for every $\epsilon>0$ there is a step function $f$ such that $\norm{c-f}_\infty<\epsilon$. As step functions have finite image, we note that $c(I)$ is paracompact. The statement follows from the fact that the closure of paracompact set in a complete metric space is compact. See also \cite[Proposition 2.16]{Pen2}. For $k>0$ it follows from continuity of $c$.

However one can give an alternate, more direct proof, which doesn't use metric structure or completeness. This approach will be useful in the next section.
Take an open cover $U_\alpha$ of the set $\overline{c(I)}$. Consider the sets $V_\alpha = \operatorname{int} c^{-1}(U_\alpha)$. Since $c$ is not continuous in general we need to apply interior operation to get open sets. However that means that open sets $V_\alpha$ might not be a cover of $I$. Indeed if $p\in I$ is not an element of $\bigcup_\alpha V_\alpha$, then $c$ is discontinuous at $p$. For each such point $p$ we add to the family $V_\alpha$ a set $]p-\epsilon,p+\epsilon[$, where $\epsilon>0$ is chosen in such a way that $c(]p-\epsilon,p+\epsilon[)$ is contained in the sum of two sets from the family $U_\alpha$. In order to find $\epsilon$ satisfying this condition take a $\alpha'$ such that $c(p)\in U_{\alpha'}$. Since $c$ is right-continuous at $p$ we can chose $\epsilon$ such that $c([p,p+\epsilon[)\subset U_{\alpha'}$. As a second step, take a set $U_{\alpha''}$ containing $\dis\lim_{x\to p^-}c(x)$. From the definition of left-side limit it follows again that we can chose $\epsilon$ such that $c(]p-\epsilon,p])\subset U_{\alpha''}$. In this way we have obtained an open cover of $I$. From compactness it follows that there exists a finite subcover of $I$. To each set in this subcover we associate one or two sets from $V_\alpha$ family which constitute a finite subcover of $\overline{c(I)}$.
\item It follows from the fact that $\mathcal{R}(I,\mathbb{E})$ is a closed subspace of the Banach space of bounded functions, cf. \cite[Proposition 2.17]{Pen2}.
\item For $k=0$ again $\mathcal{R}^0(I,\mathbb{E})$ is a closed subspace in a Banach space. For $k=1$ note that every function in $\mathcal{R}^1(I,\mathbb{E})$ is an integral of a function in $\mathcal{R}^0(I,\mathbb{E})$ (see \eqref{c^1}).
Given a Cauchy sequence of functions $c_n\in\mathcal{R}^1(I,\mathbb{E})$ we note that $c_n^{(1)}$ is a Cauchy (and thus convergent) sequence in $\mathcal{R}^0(I,\mathbb{E})$. From \cite[8.7.8]{Di} and \eqref{c^1} it follows that $c_n$ converges to integral of limit of that sequence both in norm $\norm{\;}_\infty$ and in norm $\nnorm{\;}^1_\infty$.
Thus $\mathcal{R}^1(I,\mathbb{E})$ is a Banach space and the set of primitive functions of step functions is dense in it. The proof for $k>1$ is analogous.
\item Direct consequence of definition, see also \cite[Proposition 2.18]{Pen2}.
\item See \cite[Proposition 2.21]{Pen2}.
%
%
\end{enumerate}
\endproof


In this paper we will also need the following result which is classical in the context of $C^k$-paths and which is certainly well known for specialists in the context of $k$-regulated paths. Without precise references we will give a proof of this Proposition

\begin{prop}\label{composition} Let $U$ and $V$ be two open sets in Banach spaces $\E$ and $\mathbb{F}$ respectively and $f:U\ap V$ a smooth map.
Then we have the following properties:
\begin{enumerate}
\item For any $k\geq 0$ the map $c\mapsto f\circ c$ is a smooth map from $\mathcal{R}^k(I,U)$ to $\mathcal{R}^k(I,V)$.
\item For any $k\geq 1$ the map $c\mapsto c^{(1)}$ is a smooth map from $\mathcal{R}^k(I,U)$ to $\mathcal{R}^{k-1}(I,\mathbb{E})$.
\item The map $ \mathsf{D}:c\mapsto (c(0),c^{(1)})$ is an isomorphism from $\mathcal{R}^1(I,\E)$ to $\E\times \mathcal{R}^{0}(I,\E)$.
\end{enumerate}
\end{prop}
\noindent {\bf Proof.}\\
At first if $c:I\ap U$ is a $k$-regulated path, $f\circ c: I\ap V$ is also a $k$-regulated path. 
For $k=0$ this fact is straightforward even for a map $\phi$ which is only continuous.
For $k=1$ $f\circ c$ is of class $C^0$ and as $(f\circ c)^{(1)}(t)$ one takes $f'(c(t)) c^{(1)}(t)$. The result follows from point 4 of Theorem \ref{propreg}.
For $k\geq 2$ the proof uses the same type of arguments by using an expression of the derivative of order $k$ of $f\circ c$.
%

Now consider a smooth path $\tilde{c}:\R\ap \mathcal{R}^k(I,V)$. Then we can identify $\tilde{c}$ with a map $\tilde{c}:I\times \R$ which is $C^{k-1}$ in the two variables and smooth in the second variable. Therefore the map $f\circ \tilde{c}:I\times \R\ap V$ has the same properties. But for each $\epsilon\in \R$ the map $t\mapsto f\circ \tilde{c}(\epsilon,t)$ is a $k$-regulated path in $V$ then $s\epsilon \mapsto f\circ \tilde{c}(\epsilon, \;)$ is a smooth path in $\mathcal{R}^k(I,V)$ and so $f$ is a smooth map\footnote{Since any Banach space is a convenient space a map between Banach spaces is smooth if and only if the image of a smooth path defined on $\R$ is a smooth path (cf. \cite{KrMi}).} 
which ends the proof of Point (1).

%
%

\bigskip

\noindent Point (2) follows directly from linearity of derivative $c\mapsto c^{(1)}$.

\bigskip

\noindent For Point (3), it is clear that the map $c\mapsto (c(0),c^{(1)})$ is linear and its inverse is the linear map $(x, u)\mapsto b(t)=x+\dis\int_0^t u(s)ds$. Now, we have
 
 $$\norm{c(0)}_\E+ \norm{c^{(1)}}_\infty\leq \norm{c}_\infty+\norm{c^{(1)}}_\infty \leq 2 \nnorm{c}^1_\infty$$ 
 and 
 $$\nnorm{b}^1_\infty \leq \norm{b}_\infty + \norm{b^{(1)}}_\infty\leq \norm{x}_\E + (t_1-t_0)\norm{u}_\infty + \norm{u}_\infty
 \leq (1+t_1-t_0)(\norm{x}_\E+ \norm{u}_\infty).$$
 
\noindent which shows that these linear maps are continuous and so ends the proof of Point (3).\\
\endproof

\begin{prop}\label{prop:restriction}
Let $J$ be a closed subinterval of $I$. The restriction map $\mathsf{r}_J: \mathcal{R}^k(I,\E) \ap \mathcal{R}^k(J,\E)$ is a linear contraction.
\end{prop}

\noindent {\bf Proof.}\\
Recall that the restriction of a strong regulated path has to be considered in the sense of Definition \ref{restrict}. The fact that the map $\mathsf{r}_J$ is a contraction follows directly from the definition of restriction and the norm $\nnorm{\;}^k_\infty$, see \ref{nota}:
$$\nnorm{\mathsf{r}_J(c)}^k_\infty \leq \nnorm{c}^k_\infty.$$
\endproof

\begin{prop}\label{prop:concatenation}
The concatenation map $\star: \mathcal R^0([a,b]) \times R^0([b,c]) \to R^0([a,c])$ is an isometric isomorphism of Banach spaces, where the Cartesian product is equipped with max-norm.
\end{prop}

\noindent {\bf Proof.}\\
It was previously stated that concatenation of two strong regulated paths is again a strong regulated path. The map $(\mathsf{r}_{[a,b]}, \mathsf{r}_{[b,c]})$ is its inverse. Note that even though the value of the first curve in point $b$ is seemingly lost, it reconstructed by $\mathsf{r}_{[a,b]}$ by requirement that the result be left-continuous at $b$. Directly from the definition of norm we have
$$\nnorm{c_1\star c_2}^0_\infty = \max(\nnorm{c_1}^0_\infty, \nnorm{c_2}^0_\infty).$$
\endproof

\section{$k$-regulated paths in a Banach manifold}
\label{Banachmanifold}
We now consider a Banach manifold $M$ modeled on a Banach space $\mathbb{M}$. We have the same definitions of regulated and strong regulated paths: 

\begin{defi}\label{stweak}${}$
\begin{enumerate}
\item A path $c: I\ap M$ is called regulated 
 if for any $t\in[t_0,t_1[$ the limit $\dis\lim_{s\rightarrow t^+}c(s)$ exists and for any $t\in ]t_0,t_1]$ the limit $\dis\lim_{s\rightarrow t^-_1}c(s)$ exists.
\item A regulated path $c: I\ap M$ is called strong regulated 
 if for any point $t\in [t_0,t_1[$ we have $c(t)=\dis\lim_{s\rightarrow t^+}c(s)$ and moreover $c(t_1)=\dis\lim_{s\rightarrow t_1^-}c(s)$.
\end{enumerate}
\end{defi}
 
\noindent Again for any regulated path $c:I\ap M$, there exists a unique strong regulated path $\bar{c}$ whose discontinuity set is no larger than the discontinuity set of $c$ and it coincides with $c$ outside of that set.
Just like for regulated paths in a Banach space, if $c:I\ap M$ is strong regulated, {\it the restriction of $c$} to any subinterval $J= [s_0,s_1]$ of $I$ will be strong regulated path $c_J$ defined by the restriction of $c$ to $J$. If $c:[t_0,t_1]\ap \mathbb{E}$ and $c':[t'_0,t'_1]\ap \E$ are two strong regulated paths then {\it the concatenation $c{\star}c'$ } is defined in the same way as for strong regulated paths in a Banach space and so it is also a strong regulated path.

\bigskip

A partition of $I=[t_0,t_1]$ is a non-decreasing sequence $\t=(\t_i)_{i=0,\dots,n}$ such that $t_0=\t_0\leq \t_1\leq\cdots\leq\t_i\leq\cdots\leq \t_{n}=t_1$.
Of course if $c$ is strong regulated then for any partition $\t=(\t_i)_{i=0,\dots,n}$ of $[t_0,t_1]$ the restriction $c_i=c_{[\t_{i-1},\t_{i}]}$ is strong regulated for all $i=1,\dots,n$.

\bigskip

{\bf From now on in this section we fix an atlas $\mathfrak{A}=(U_\alpha, \phi_\alpha)_{\alpha \in A}$ on $M$ such that each domain $U_\alpha$ is connected and simply connected.}

\bigskip

\begin{defi}\label{kregM} Let $c: I\ap M$ be a path.
A $0$-regulated path is a strong regulated path. For $k>0$ $c$ is a $k$-regulated path if and only if
for any chart $(U_\alpha,\phi_\alpha)$ such that $U_\alpha\cap c(I)\not=\emptyset$ and for any subinterval $J\subset I$ such that $c(J)\subset U_\alpha$, the path $\phi_\alpha\circ c_{J}$ is a $k$-regulated path in $\mathbb{M}$.
 \end{defi}

 \begin{rem}\label{inddefi-smooth}${}$\normalfont 
 \begin{enumerate}
 \item The definition of $k$-regulated path does not depend on the choice of the atlas $\mathfrak{A}$.
 \item If $c: I\ap M$ is $k$-regulated for all $k\in \N$, it is smooth in the usual sense.
 \end{enumerate}
 \end{rem}
 
\begin{prop}\label{repar} Let $c:I\ap M$ be a $k$-regulated path.
\begin{enumerate}
\item The closure of $c(I)$ is a compact subset of $M$.
\item Let $c:I\ap M$ be a $k$-regulated path. For any increasing map $h:[s_0,s_1]\ap I$ of class $C^{k+1}$ the path $c\circ h:[s_0, s_1]\ap M$ is a $k$-regulated path called {\it a reparametrization of $c$}. In particular for any $k$-regulated path $c:I\ap M$, there exists a canonical increasing affine map $h:[0,1]\ap I$ such that $c\circ h:[0,1]\ap M$ is a $k$-regulated path.
\end{enumerate}
\end{prop}

\noindent {\bf Proof.}${}$
\begin{enumerate}
\item Note that for $k\geq 1$ the path $c$ is continuous and in consequence it is obvious. For $k=0$ the alternate proof from Theorem \ref{propreg} works in this case without any modifications.


\item Note that the first part of Point 2 is true for $M=\M$ according to Theorem \ref{composition} Point (1). Therefore the definition of a $k$-regulated path in a manifold implies clearly this result. Now the map 
$$s\mapsto {s}(t_1-t_0)+t_0$$
is smooth strictly increasing from $[0,1]$ onto $I=[t_0,t_1]$. The last property is a direct consequence of the first part.\\
\end{enumerate}

\endproof

\section{$0$-regulated lifts of $1$-regulated paths }


Let $\pi:E\ap M$ a Banach bundle on a Banach manifold modeled on a Banach space $\M$ with fibers modeled on a Banach space $\E$. We denote by $p_E:TE\ap E$ (resp. $p_M:TM\ap M$) the tangent bundle of $E$ (resp. $M$).

\begin{defi}\label{bpath} Let $\g:I\ap M$ be a path.
A $k$-regulated path $c:I\ap E$ is called a $k$-regulated lift or a $k$-regulated section along $\g$ if we have $\pi\circ c=\g$. If $E=TM$, a $0$-regulated section over $\g$ will be called a $k$-regulated vector field along $\g$.
\end{defi}

According to Remark \ref{repar} any $0$-regulated lift $c:I\ap E$ of a $1$-regulated path $\g:I\ap M$, can be reparametrized as a $0$-regulated lift $\bar{c}:[0,1]\ap E$ over a $1$-regulated path $\g\circ h:[0,1]\ap M$. Therefore we can assume that $1$-regulated paths $\g$ are defined on $[0,1]$ and so all their $0$-regulated lifts $c$ will also be defined on $[0,1]$.

\begin{nota} \label{VFg}\normalfont {\bf --- from now $I$ will denote the interval $[0,1]$}.
\end{nota}

The set of $1$-regulated paths in $M$ defined on $I$ is denoted $\mathsf{P}(M)$ and the set of $0$-regulated lifts (to $E$) of paths in $\mathsf{P}(M)$ is denoted by $\mathcal{P}(E)$.
Of course $\pi:E\to M$ induces a surjective map (also denoted by $\pi$) from $\mathcal{P}(E)$ onto $\mathsf{P}(M)$.

\subsection{Banach structure on the set of $0$-regulated sections and $0$-regulated vector fields along $1$-regulated paths}\label{GcGg}

We denote by $\G(\g)$ the set of $1$-regulated vector fields along $\g\in \mathsf{P}(M)$ and by $\G(c)$ the set of $0$-regulated vector fields along $c\in \mathcal{P}(E)$.
 

\bigskip

Our purpose is to give some nice parametrization of $\G(\g)$ and $\G(c)$ which provides these sets with a Banach structure. For this objective we need the following technical result:

\begin{lem}\label{parametre}
Consider $c\in \mathcal{P}(E)$ and $\g=\pi\circ c$.
Let $\g^*(E)$ be the pull-back of $E$ over $\g$.
\begin{enumerate}
\item There exists a homeomorphism $G_\g:I\times \E\ap \g^*(E)$ which is a bundle isomorphism over the $Id$ of $I$, which is $1$-regulated in the first variable and smooth in the other variable such that the following diagram is commutative
 
\begin{equation}\label{Gg}
\begin{tikzcd}
I\times \E\arrow{r}{G_\g}\arrow{d}{}
&\g^*(E)\arrow{r}{\bar{\g}}\arrow{d}{{\pi_\g}}
&E\arrow{d}{\pi}\\
I\arrow{r}{Id}&I\arrow{r}{\g}&M
\end{tikzcd}
\end{equation}
\item Let $\bar{G}_\g= \bar{\g}\circ G_\g$ and consider the pull-back $(\bar{G}_\g)^*(TE)\ap I\times \E$ of the bundle $p_E:TE\ap E$. Then there exists a homeomorphism $H_\g: I\times\E\times\M\times \E\ap \bar{G}_\g^*(TE)$, which is a bundle isomorphism over $Id$ of $I\times\E$ $1$-regulated in the first variable and smooth in the other variables, and such that the following diagram is commutative


\begin{equation}\label{hatGg}
\begin{tikzcd}
I\times \E\times\M\times\E \arrow{r}{H_\g}\arrow{d}{q_{I\times\E}}
&\bar{G}_\g^*(TE)\arrow{r}{\widehat{G_\g}}\arrow{d}{{p_{G_\g}}}
 & TE\arrow{d}{p_E}\\
I\times\E\arrow{r}{Id}\arrow{d}&I\times\E\arrow{r}{\bar{G}_\g}\arrow{d}& E\arrow{d}{\pi}\\
I\arrow{r}{Id}&I\arrow{r}{\g}&M\\
\end{tikzcd}
\end{equation}
\end{enumerate}
\end{lem}

\noindent{\bf Proof}${}$\\
(1) Since $I$ is (smooth) contractible and $\g=\pi\circ c$ is a $1$-regulated path, there exists a bundle isomorphism $G_{{\g}}$ of class $C^{0}$ over the identity from the trivial bundle 
$ I\times\E\ap I$ to the bundle ${\g}^*(E)\ap I $ ({cf.} \cite{AMR} Theorem 3.4.35). In fact since $\g$ is $1$-regulated, it follows that $G_\g$ is also $1$-regulated in the variable $t\in I$. 
Smoothness in the other argument follows from linearity in the fibers.
Since $G_\g$ is a bundle isomorphism the commutativity of the diagram follows.

\noindent (2) 
By the same argument as previously, there exists a bundle isomorphism $H_\g$ from the trivial bundle $I\times \E\times\M\times\E\ap I\times \E$ to the bundle $(\bar{G}_\g)^*(TE)\ap I\times \E$. The fact that $\bar{G}_\g$ is $1$-regulated in the variable $ t$ and smooth in the other variables implies the corresponding properties of $H_\g$. Now since $H_\g$ is a bundle isomorphism, the commutativity of diagrams follows.\\
\endproof

 If $q_\E$ is the projection of $I\times\E$ on $\E$, the map $q_\E\circ G_\g^{-1}: \g^*(E)\ap \E$ is called a {\it representation of $\g^*(E)$}. In particular for any $0$-regulated lift $c:I\ap E$ of $\g$, then 
${u}=q_\E\circ G_g^{-1}(c)$ belongs to $\mathcal{R}^0(I,\E)$. The pair $(\g, {u})$ will be called a {\it representation of $c$}. According to the notations of Lemma \ref{parametre}, a map $\vartheta: I\ap E$ induces a section of $\g^*(E)$ if an only if $\pi\circ \vartheta=\g$. Then, from
Lemma \ref{parametre} Point (1), $\mathcal{S}_\g=q_\E\circ \bar{G}_\g^{-1}$ induces a vector isomorphism from the vector space $\G(\g^*(E))$ of $0$-regulated lifts $c:I\ap E$ along $\g$ onto the Banach vector space $\mathcal{R}^0(I,\E)$. \\

By the same arguments applied to $E=TM$ we also have a vector isomorphism $\mathcal{T}_\g$ from the vector space $\G(\g)$ of $1$-regulated vector fields along $\g$ to $\mathcal{R}^1(I,\M)$.\\

Using the bundle morphism $\widehat{G_\g}: \bar{G}_\g^*(TE)\ap TE$ over $\bar{G}_\g$ from the diagram \eqref{hatGg} and considering a $0$-regulated lift $c: I\ap E$ of $\g$ 
we note that any map $\mathcal{X}:I\ap TE$ satisfies $p_E\circ \mathcal{X}=c$ if and only if 
$q_{I\times \E}\circ H_\g^{-1}\circ \widehat{G_\g}^{-1}\circ \mathcal{X}=\bar{G}_\g^{-1}\circ c$.
 
 Given a vector field $\mathcal{X}$ along $c$, if $q_{\M\times\E}$ is the projection of $I\times\E\times\M\times\E$ on $\M\times\E$, then the map 
 $$t\mapsto q_{\M\times \E}\circ H_\g^{-1}\circ\widehat{G_\g}^{-1}\circ \mathcal{X}(t)$$
 is map from $I$ to $\M\times \E$ which is a $0$-regulated path if and only if $\mathcal{X}$ is a $0$-regulated path. Note that $T\pi (\mathcal{X})$ is vector field along $\g$ and according to the commutativity of the diagrams in Lemma \ref{parametre}, the projection of 
$$t\mapsto q_{\M\times \E}\circ H_\g^{-1}\circ \widehat{G_\g}^{-1}\circ \mathcal{X}(t)$$
 on $\M$ is nothing else but $t\mapsto\mathcal{T}_\g(T_{c(t)}\pi(\mathcal{X}(t))$.\\

Let $\mathcal{T}_c$ be the map from the vector space $\G(c)$ to $\mathcal{R}^1(I,\M)\times \mathcal{R}^0(I,\E)$ defined by 
$$\mathcal{X}\mapsto q_{\M\times\E}\circ H_\g^{-1}\circ \widehat{G_\g}^{-1}(\mathcal{X}).$$
 Clearly $\mathcal{T}_c$ is a vector isomorphism and if $q_1$ is the projection of $\mathcal{R}^1(I,\M)\times \mathcal{R}^0(I,\E)$ onto $\mathcal{R}^1(I,\M)$, we have
$q_1\circ \mathcal{T}_c(\mathcal{X})=\mathcal{T}_\g(T\pi\circ \mathcal{X}).$ Therefore any vector field $\mathcal{X}$ along $c$ will be written as a pair $(\varphi,\vartheta)=\mathcal{T}_c(\mathcal{X})$ with $\varphi\in \mathcal{R}^1(I,\M)$ and $\vartheta\in \mathcal{R}^0(I,\E)$. Note that we have of course $\varphi=\mathcal{T}_\g(T\pi(\mathcal{X}))$.
Therefore {\it $(\varphi,\vartheta)$ will be called a representation of $\mathcal{X}$}. Thus for one choice of the trivializations ${G}_\g$ and ${H}_\g$, we get the following commutative diagram
\begin{equation}\label{locrepfib}
\begin{tikzcd}
\G(\g^*(E))\arrow{r}{\mathcal{S}_\g}&\mathcal{R}^0(I,\E)\\
\G(c)\arrow{r}{\mathcal{T}_c}\arrow{u}{p_E}\arrow{d}{T\pi}
&\mathcal{R}^1(I,\M)\times \mathcal{R}^0(I,\E)\arrow{u}{q_2}\arrow{d}{q_1}\\
\G(\g)\arrow{r}{\mathcal{T}_\g}&\mathcal{R}^1(I,\M)
\end{tikzcd}
\end{equation}

Representations $\mathcal{S}_\g$, $\mathcal{T}_\g$ and $\mathcal{T}_c$ of $\G(\g^*(E))$, $\G(\g)$ and $\G(c)$ such that the diagram \eqref{locrepfib} is commutative 
are called {\it compatible}. Note that if $\mathcal{S}'_g$, $\mathcal{T}'_\g$ and $\mathcal{T}'_c$ are other compatible representations of $\G(\g^*(E))$, $\G(\g)$ and $\G(c)$ then there exists a $0$-regulated field $A^E_\g$, a $1$-regulated field $A_{\g}$ and a $0$-regulated field $A_{c}$ of automorphisms of $\E$, $\M$ and $\M\times\E$ respectively such that: 
\begin{equation}\label{compat}
 \mathcal{S}'_{\g}=A^E_{\g}\circ \mathcal{S}_{\g},\;\; \mathcal{T}'_{\g}=A_{\g}\circ \mathcal{T}_{\g},\;\; \mathcal{T}'_{c}=A_{c}\circ \mathcal{T}_{c} \textrm{ and } A_{\g}=q_1\circ A_{c}.
\end{equation}

\begin{nota}\label{triviaVF} ${}$\\ \normalfont Fixing a choice of such compatible trivializations $\mathcal{T}_\g$ and $\mathcal{T}_c$ of $\G(\g)$ and $\G(c)$ we will identify $\mathcal{X}\in \G(c)$ with the pair $(\varphi,\vartheta)=\mathcal{T}_c(\mathcal{X})$ and $X=T\pi(\mathcal{X})$ with $\varphi$. In this way $\vartheta$ can be considered as a $0$-regulated vertical vector field of $E$ over $c$ as well as a $0$-regulated section of $E$ over $\g$ according to the fact that a vertical vectors at $(x,u)\in E$ can be canonically
identified with vectors of the fiber $E_x$.
\end{nota}

\begin{prop}\label{Banachstr} With the previous notations the representations $\mathcal{S}_\g$, $\mathcal{T}_\g$ and $\mathcal{T}_c$ induce Banach structures
$\norm{\;}_{\mathcal{S}_\g}$, $\norm{\;}_{\mathcal{T}_\g}$, $\norm{\;}_{\mathcal{T}_c}$ on $\G(\g^*(E))$, $\G(\g)$ and $\G(c)$ from $\mathcal{R}^0(I,\E)$, $ \mathcal{R}^1(I,\E)$ and $\mathcal{R}^1(I,\M)\times \mathcal{R}^0(I,\E)$ which are independent on the choice of compatible trivializations of $\g^*(E)$, $\g^*(TM)$ and $\bar{G}_\g^*( TE)$.
\end{prop}

 \begin{rem}\label{localtrivial}${}$
\begin{enumerate}
 \normalfont
\item Assume that for a $1$-regulated path $\g:I\ap M$ there exists a chart $(U,\phi)\in \mathfrak A$ such that $\g(I)$ is contained in $U$ and that $E_{|U}$ is trivializable and let $\Phi:E_{| U}\ap \phi(U)\times \E$ be a trivialization. Then $\Phi$ induces a trivialization $G_\g: I\times \E \to \g^*(E)$ of the bundle $\g^*(E)$ by $G_\g(t,u)=\Phi^{-1}(\phi(\g(t)), u)$. Then for any section $c$ of $\g^*(E)$, $\Phi(c)$ is a pair of paths $(x,u): I\ap \M\times \E$ and the representation $\mathcal{S}^\phi_\g$ of $\G(\g^*(E))$ into $\mathcal{R}^0(I,\E)$ associated to $G_\g$ is exactly the projection on $\E$ of $\Phi(c)$. In the same way a natural representation $\mathcal{T}^\phi_\g$ of $\G(\g)$ into $\mathcal{R}^1(I,\M)$ is the projection on $\M$ of $T\phi$. Finally, we also have a natural representation $\mathcal{T}^\phi_c$ of $\G(c)$ into $\mathcal{R}^1(I,\M)\times \mathcal{R}^0(I,\M)$ which is a projection on $\M\times \E$ of $T\Phi$.
 
\item Let $\t=(\t_i)_{i=0,\dots,n}$  be a partition of $I$, 
$(U_i,\phi_i)_{i=1,\dots,n}\subset \mathfrak A$ be a covering of a $1$-regulated path $\g:I\ap \M$ 
such that $\overline{\g([\t_{i-1}, \t_{i}])}\subset U_i$ 
and $c:I\ap E$ be a $0$-regulated lift of $\g$. We set $\g_i=\g_{|[\t_{i-1},\t_i]}$ and $c_i=c_{|[\t_{i-1},\t_i]}$. Then a trivialization $G_\g:I\times \M\ap \g^*(E)$ (resp. $\bar{G}_\g: I\times \M\times\E\ap c^*(TE)$) induces by restriction a trivialization $G_{\g_i}:[\t_{i-1},\t_i]\times \M\ap \g_i^*(E)$ (resp. $\bar{G}_{\g_i}: [\t_{i-1},\t_i]\times \M\times\E\ap c_i^*(TE)$). This implies that the corresponding trivialization $\mathcal{T}_{\g_i}:\G(\g_i)\ap \mathcal{R}^1([\t_{i-1},\t_i], \M)$ (resp. $\mathcal{T}_{c_i}:\G(c_i)\ap \mathcal{R}^1([\t_{i-1},\t_i], \M)\times \mathcal{R}^0([\t_{i-1},\t_i], \E)$) is the restriction of $\mathcal{T}_\g$ (resp. $\mathcal{T}_c$) to $\G(\g_i)$ (resp. $\G(c_i)$). Moreover we have a $1$-regulated (resp. $0$-regulated) field $t\mapsto A_{\g_i}$ (resp. $t\mapsto A_{c_i}$) of automorphisms of $\M$ (resp. $\M\times\E$ ) such that ({cf.} (\ref{compat})):
$$\mathcal{T}_{\g_i}^{\phi_i}=A_{\g_i}\circ \mathcal{T}_{\g_i} ,\;\; \mathcal{T}_{c_i}^{\phi_i}=A_{c_i}\circ \mathcal{T}_{c_i} \textrm{ and } A_{\g_i}=q_1\circ A_{c_i}.$$

\item Given a trivialization $G_\g :I\times\M\ap \g^*(TM)$, obtained from Lemma \ref{parametre} in case $E=TM$, for each $s\in I$ the linear map ${T}_s:\M\equiv T_{\g(0)}M\ap T_{\g(s)}M$ defined by $T_s(v)=\bar\g\circ G_\g(s,v)$ is an isomorphism of Banach spaces. Set $P_t^s=T_t\circ T_s^{-1}$. 
 
Given any $X\in \G(\g)$ we can define 
$$\nabla_{\g}X(t)=\dis\frac{d}{dt}\{P_s^t(X(s))\}_{| s=t}$$
outside at most countable subset of $I$. Using local chart and Point 1, it is easy to see that it can be extended to a $0$-regulated vector field along $\g$. It follows that the norm $\norm{\;}_{\mathcal{T}_\g}$ is equivalent to the norm 
\begin{equation}\label{nablaN}
\norm{X(0)}_\M+\norm{\nabla_\g X}_\infty.
\end{equation}
Operators $P^s_t$ play the role of a parallel transport. Namely if we would have a linear connection on $M$ and $\nabla$ would be its Koszul operator, then we would have a parallel transport along $\g$ defined in the same way as previously and (\ref{nablaN}) would be the standard norm on $\G(\g)$ associated to $\nabla$ (see for instance \cite[A8]{Sch}). 
\end{enumerate}
\end{rem}

\noindent{\bf Proof of Proposition \ref{Banachstr}}${}$\\
\noindent 
Since $\mathcal{T}_\g$ is a vector bundle isomorphism from $\G(\g)$ to $\mathcal{R}^1(I,\M)$ we can provide $\G(\g)$ with the pull-back Banach structure. Then $\mathcal{T}_\g$ induces an isometry between these two Banach spaces. Indeed if $\varphi=\mathcal{T}_\g(X)$ then outside at most a countable subset of $I$ we have 
$$\dot{\varphi}(t)=\dis\frac{d}{dt}\{\mathcal{T}_\g\}X(t)+\mathcal{T}_\g\dis\frac{d}{dt}\{X\}(t).$$
Since $\varphi$ is $1$-regulated on $\G(\g)$ the induced norm is
\begin{equation}\label{inducedN}
\norm{X}_{\mathcal{T}_\g}=\max\left\{\norm{\mathcal{T}_\g(X)}_\infty, \norm{\dis\frac{d}{dt}\{\mathcal{T}_g\}X(t)+\mathcal{T}_g\dis\frac{d}{dt}\{X\}(t)}_\infty\right\}.
\end{equation}

From Lemma \ref{parametre}, if $G'_\g:I\times \E\ap \g^*(TM)$ is another trivialization, then $G'_\g\circ G_\g^{-1}$ is an automorphism of the bundle $I\times \M\ap I$. As we have already seen, there exists $1$-regulated field $t\mapsto A(t)$ from $I$ to $GL(\M)$ such that $\mathcal{T}_\g(X(t))=A(t)\mathcal{T}'_\g(X(t))$ for all $t\in I$. Since any $1$-regulated path is bounded (cf. Proposition \ref{repar}), if $\varphi=\mathcal{T}_\g(X)$ and $\varphi'=\mathcal{T}'_\g(X)$ then it is easy to show that there exists $K>0$ such that we have

$$\norm{\varphi}_\infty\leq {K}\norm{\varphi'}_\infty \qquad \textrm{and} \qquad \norm{\dot{\varphi}'}_\infty\leq K \sup\{\norm{\varphi'}_\infty,\norm{\dot{\varphi}'}_\infty\}.$$

\noindent This implies that the norms on $\G(\g)$ induced by $\mathcal{T}_\g$ and $\mathcal{T}'_\g$ are equivalent. Since any $0$-regulated path is also bounded, we can apply formally the same arguments to the isomorphism $\mathcal{S}_\g$ from $\G(\g^*(E))$ to $\mathcal{R}^0(I,\E)$.

Finally, for the vector space $\G(c)$ again we can use analogous arguments as the previous ones for two trivializations $H_\g$ and $H'_\g$ of the bundle $\bar{G}_\g^*(TE)$ over $\bar{G}_\g^{-1}(c)$. 

\endproof

\subsection{Topology on the set of regulated lifts of paths}\label{topologyPE}
Since any $0$-regulated path defined on $I$ has a compact range, we can provide the set $\mathcal{P}(E)$ with the topology generated by sets $\mathcal{N}(K,U,V,W)$ of paths $ c \in \mathcal{P}(E) $ such that the closure of $c(K)$ is contained in $ U$, closure of $(\pi \circ c)(K)$ is contained in $V$ and closure of $(\pi \circ c)^{(1)}(K)$ is contained in $W$, where $K$ is a compact subset of $I$, $U$ is an open set in $E$, $V$ is an open set in $M$ and $W$ is an open set in $TM$.

{\it This topology is Hausdorff}. Indeed, if $c_1\not=c_2$, there exists $\theta\in I$ such that $c_1(\theta)\not=c_2(\theta)$. Now take $K=\{\theta\}$. Since $E$ is Hausdorff, the points $c_1(\theta)$ and $c_2(\theta)$ have open disjoint neighbourhoods $U_1$ and $U_2$. Sets $\mathcal{N}(K,U_1,M,TM)$ and $\mathcal{N}(K,U_2,M,TM)$ are required open disjoint neighbourhoods of $c_1$ and $c_2$.

We can provide $\mathsf{P}(M)$ with the topology in the same way that is generated by sets of type $\mathcal{N}(K,V,W)$ of paths $\g \in \mathsf{P}(M)$ such that the closure of $\g(K)$ is contained in $V$ and the closure of $\g^{(1)}(K)$ is contained in $W$ where $K$ is a compact subset of $I$ and $V$ is an open set in $M$ and $W$ an open set in $TM$. Now the map $\pi:E\ap M$ induces a natural map $c\mapsto \pi\circ c$ from $\mathcal{P}(E)$ to $\mathsf{P}(M)$ which we denote also $\pi$. It is continuous and open with respect to these topologies. This implies in particular that $\mathsf{P}(M)$ is a Hausdorff space.

\bigskip
 
If $(U,\phi)$ is a chart domain such that there exists a trivialization $\Phi:E_{|U}\ap \phi(U)\times \E$, then $\mathcal{P}(E_{| U})$ is an open set of $\mathcal{P}(E)$. Now we have:
 
\begin{prop}\label{basecont}${}$
\begin{enumerate}
\item For $i=1,2$ let $\pi_i:E_i\ap M_i$ be two Banach bundles and $F$ be a bundle morphism from $E_1\ap M_1$ to $E_2\ap M_2$ over $f:M_1\ap M_2$. Then the map $F_*:\mathcal{P}(E_1)\ap \mathcal{P}(E_2)$ (resp. $f_*:\mathcal{P}(M_1)\ap \mathcal{P}(M_2)$) defined by $F_*(c)=F\circ c$ (resp. $f_*(\g)=f\circ \g$ ) is continuous. Moreover $F_*$ and $f_*$ are injective (resp. a homeomorphism) if $F$ is injective (resp. isomorphism). Moreover we have $\pi_2\circ F_*=f_*\circ\pi_1$.
\item The map $\Phi_*: \mathcal{P}(E_{| U})\ap \mathcal{R}^1(I,\phi(U))\times \mathcal{R}^0(I,\E)$ is a homeomorphism.
\item The evaluation map ${\bf Ev^t}:\mathcal{P}(E)\ap E$ defined by ${\bf Ev^t}(c)=c(t)$ is continuous.
\end{enumerate}
\end{prop}

\noindent {\bf Proof}${}$\\
(1)${}$ For $c\in \mathcal{P}(E_1)$ take an open neighbourhood $\mathcal{N}_2(K,V,W)$ of $F\circ c$ in $\mathcal{P}(E_2)$. Consider $\mathcal{N}_1(K, F^{-1}(V), Df^{-1}(W))$. Since $F$ and $f$ are smooth, it follows that $F^{-1}(V)$ and $Df^{-1}(W)$ are open sets of $E_1$ and $TM_1$ respectively. Therefore $\mathcal{N}_1(K, F^{-1}(V), Df^{-1}(W))$ is a neighbourhood of $c$ and $F_*( \mathcal{N}_1(K, F^{-1}(V), Df^{-1}(W))$ is contained in $\mathcal{N}_2(K,V,W)$. Thus $F_*$ is continuous. It is clear that if $F$ is injective so is $F_*$. Now if $F$ is an isomorphism, by application of the previous result to the isomorphism $F^{-1}$ we obtain that $F_*$ is a homeomorphism. All the results about $f_*$ comes from the relation $\pi_2\circ F=f\circ \pi_1$.\\
 
\noindent (2)${}$ Note that the set $\mathcal{P}(\M\times\E)$ is exactly $\mathcal{R}^1(I,\M)\times\mathcal{R}^0(I,\E)$. Moreover since $\M$ and $\E$ are metrizable, the topology defined on $\mathcal{P}(\M\times\E)$ coincides with product topology of the Banach topology of $\mathcal{R}^1(I,\M)$ and $\mathcal{R}^0(I,\E)$. 
Now $\phi(U)$ is an open set of $\M$ so $\mathcal{R}^1(I,\phi(U))\times \mathcal{R}^0(I,\E)$ is an open set of $\mathcal{R}^1(I,\M)\times \mathcal{R}^0(I,\E)$. Thus the proof of (2) is a consequence of (1) by composition of homeomorphisms.\\
 
\noindent (3)${}$ Take any open set $V$ in $E$. Then $({\bf Ev^t})^{-1}(V)$ is the set $\{c\in \mathcal{P}(E)\;|\; c(t)\in V\}$. If we take the compact interval $K=\{t\}$, it follows that by construction of the topology on $\mathcal{P}(E)$ that $({\bf Ev^t})^{-1}(V)$ is an open set.\\
\endproof

 

\section{Banach manifold structure of $\mathcal{P}(E)$ }\label{PEn}

\emph{In this section, we fix an atlas $\mathfrak{A}=\{U_\alpha,\phi_\alpha\}_{\alpha\in A}$ on $M$ such that each $U_\alpha$ is contractible. In this situation recall that from Theorem 3.4.35 in \cite{AMR} it follows that the bundle $E_{|U_i}$ is trivializable}.\\

We have the following basic fundamental result:

\begin{theo} \label{banachP}${}$
\begin{enumerate}
\item The set $\mathsf{P}(M)$ has a Banach manifold structure modeled on the Banach space $\mathcal{R}^1(I,\M)$ such that each chart domain of this structure is open for its compact open topology. Moreover for any $\g\in \mathsf{P}(M)$ we have a Banach isomorphism from $T_\g\mathsf{P}(M)$ onto $\G(\g)$.
\item The set ${\mathcal{P}}(E)$ has a Banach manifold structure modeled on the Banach space $\mathcal{R}^1(I,\M)\times\mathcal{R}^0(I ,\E)$ such that each chart domain is open for the natural topology of ${\mathcal{P}}(E)$ and we have a Banach isomorphism from $T_c {\mathcal{P}}(E)$ onto $\G(c)$. Moreover $\pi:\mathcal{P}(E)\ap \mathsf{P}(M)$ is a Banach bundle with typical fiber $\mathcal{R}^0(I ,\E)$.
\end{enumerate}
\end{theo}
 
Before proving the theorem we begin by the following Lemma

\begin{lem}\label{isoTCGc}
Given a partition $\t=(\t_i)_{i=0,\dots,n}$ of the interval $I$ we define
the map $\mathcal{I}_\t$ which to 
$$(x, y_1, v_1,\dots,y_n,v_n)\in\M\times\prod_{i=1}^n(\mathcal{R}^0\big([\t_{i-1},\t_i],\M)\times\mathcal{R}^0([\t_{i-1},\t_i],\E)\big)$$ 
associates a pair of functions $(\varphi,\vartheta)$ in $ \mathcal{R}^{1}(I,\M)\times \mathcal{R}^0(I,\E)$ in the following way:
\begin{align*}
\vartheta &= v_1\star v_2 \star\ldots\star v_n,\\
\varphi &= x+\dis\int_0^t (y_1\star y_2 \star\ldots\star y_n)(s)ds,
\end{align*}
where $\star$ is a concatenation of curves defined in Section \ref{regulated}. We define also a map $\mathsf{I}_\t:\M\times\prod\limits_{i=1}^n\mathcal{R}^0([\t_{i-1},\t_i],\M)\to \mathcal{R}^{1}(I,\M)$ as follows
$$\mathsf{I}_\t(x, v_1,\dots,v_n) = \varphi.$$
Then the maps $\mathcal{I}_\t$ and $\mathsf{I}_\t$ are Banach isomorphisms.
\end{lem}

The lemma follows by applying $n$ times Proposition \ref{prop:restriction} and Proposition \ref{composition} Point 3.

\medskip


Consider a curve $c\in \mathcal{P}(E)$, i.e. $0$-regulated path $c:I\ap E$ such that $\g=\pi\circ c$ is $1$-regulated. Assume for some subinterval $J=[s_0,s_1]$ of $I$ that $\g(J)$ is contained in the domain $U$ of a chart $(U,\phi)$ such that $E_{| U}$ is trivializable by a map $\Phi: E_{| U}\ap \phi(U)\times \E\subset \M\times \E$. Then the image $\Phi(c_{| J})$ can be written as a pair $ (\phi\circ \g_{| J},u): J\ap \M\times \E$ and since $\g$ is $1$-regulated, we have ({cf.} (\ref{c^1})):
$$\phi( \g_{| J}(t))=\phi(\g(s_0))+\dis \int_{s_0}^t (\phi\circ \g)^{(1)}(s)ds.$$
In such a situation, we set 
\begin{equation}\label{n=1}
\tilde{\phi}(\g_{|J})=(\phi\circ \g_{| J})^{(1)}\;\;\;\tilde{\Phi}(c_{| J})=(\tilde{\phi}(\g),u ).
\end{equation}
Therefore the path $\tilde{\phi}(\g_{|J})$ belongs to $\mathcal{R}^0([s_0,s_1], \M)$ and the path $\tilde{\Phi}(c_{| J})$ belongs to $\mathcal{R}^0([s_0,s_1], \M)\times\mathcal{R}([s_0,s_1], \E)$ and of course 
\begin{equation}\label{n=1proj}
\tilde{\phi}(\pi\circ c_{| J})=p_\M\circ \tilde{\Phi}(c_{| J}).
\end{equation}

{\bf We come back to the general context of the proof of Theorem \ref{banachP}. Recall that the atlas $\mathfrak{A}=\{U_\alpha,\phi_\alpha\}_{\alpha\in A}$ of $M$ is fixed and each $U_\a$ is assumed to be connected and simply connected and so there exist a trivialization $$\Phi_\alpha: E_{U_\alpha}:=E_{| U_\alpha}\to \phi_\alpha(U_\alpha)\times \mathbb{E}\subset \mathbb{M}\times \mathbb{E}.$$}
For each integer $n\geq 1$, each partition $\t=(\t_i)_{i=0,\dots,n}$ and each family $\mathfrak{U}=(U_i)_{i=1,\dots,n}$ of chart domains of $\mathfrak{A}$ we consider the sets:

$$\mathsf{P}(\t,\mathfrak{U})=\left\{\g\in {\mathsf{P}}(M) \;|\; {\g([\t_{i-1},\t_i])}\subset U_i, \forall \; i=1,\dots n \right \},$$
$$\mathcal{P}(\t,\mathfrak{U})=\left\{c\in {\mathcal{P}}(E) \;|\; \pi\circ c \in \mathsf{P}(\t,\mathfrak{U})\right \}.$$

Since all functions in $\mathsf{P}(M)$ are continuous and thus have compact image,
the set of all such ${\mathcal{P}}(\t,\mathfrak{U})$ (resp. $\mathsf{P}(\t,\mathfrak{U})$) for all integers $n$, all partitions $\t=(\t_i)_{i=0,\dots, n}$ and all sequences $\mathfrak{U}=(U_i)_{i=1,\dots,n}$ is a covering of $\mathcal{P}(E)$ (resp. $\mathsf{P}(M)$).\\

\noindent According to (\ref{n=1}), for any set ${\mathcal{P}}(\t,\mathfrak{U})$ and $\mathsf{P}(\t,\mathfrak{U})$ we consider the maps
\begin{equation}\label{Phin}
\mathcal{F}_{\t,\mathfrak{U}}(c)=(\phi_1(\pi\circ c(0)),\tilde{\Phi}_1(c_1),\dots,\tilde{\Phi}_{n}(c_{n})),
\end{equation}
where $c\in {\mathcal{P}}(\t,\mathfrak{U})$, $\g=\pi\circ c$,$\;c_i=c_{[\t_{i-1},\t_{i}]}$, and

\begin{equation}\label{Phin2}
\mathsf{f}_{\t,\mathfrak{U}}(\g)=(\phi_1(\g(0)),\tilde{\phi}_1(\g_1),\dots,\tilde{\phi}_{n}(\g_{n})),
\end{equation}
where $\g\in\mathsf{P}(\t,\mathfrak{U})$, $\g_i=\g_{[\t_{i-1},\t_{i}]}$.\\

Therefore $\mathcal{F}_{\t,\mathfrak{U}}(c)$ belongs to $\M\times \dis\Pi_{i=1}^{n}\mathcal{R}^{0}([\t_{i-1},\t_{i}] ,\M)\times\mathcal{R}^0([\t_{i-1},\t_{i}] ,\E)$ and $\mathsf{f}_{\t,\mathfrak{U}}(\g)$ belongs to $\M\times \dis\Pi_{i=1}^{n}\mathcal{R}^{0}([\t_{i-1},\t_{i}] ,\M)$. \\

\begin{com}\label{atlassimple}\normalfont
According to Lemma \ref{isoTCGc}, first we must show that the set of pairs $\left(\mathcal{P}(\t,\mathfrak{U}), \mathcal{I}_\t\circ \mathcal{F}_{\t,\mathfrak{U}}\right)$ 
(resp. $\left(\mathsf{P}(\t,\mathfrak{U}), \mathcal{I}_\t\circ \mathsf{f}_{\t,\mathfrak{U}}\right)$) for all integers $n$, all partitions $\t=(\t_i)_{i=0,\dots, n}$ and all sequences $\mathfrak{U}=(U_i)_{i=1,\dots,n}$, defines a Banach manifold structure on $\mathcal{P}(E)$ (resp. $\mathsf{P}(M)$).

However since $\mathcal{I}_\t$ and $\mathsf{I}_\t$ are isomorphisms, it is sufficient to show that each $c\in {\mathcal{P}}(E)$ belongs to a chart domain $\mathsf{P}(\t,
\mathfrak{U})$ which is homeomorphic to an open set in
$\M\times \dis\Pi_{i=1}^{n}\mathcal{R}^{0}([\t_{i-1},\t_{i}] ,\M)\times\mathcal{R}([\t_{i-1},\t_{i}] ,\E)$ via $\mathcal{F}_{\t,\mathfrak{U}}$ and the transition maps between two such 
chart $\left(\mathcal{P}(\t,\mathfrak{U}), \mathcal{F}_{\t,\mathfrak{U}}\right)$ and $\left(\mathcal{P}(\s,\mathfrak{V}), \mathcal{F}_{\s,\mathfrak{V}}\right)$ are local 
diffeomorphisms between two Banach spaces of type $\M\times \dis\Pi_{i=1}^{m}\mathcal{R}^{0}([\s_{j-1},\s_{j}] ,\M)\times\mathcal{R}^0([\s_{j-1},\s_{j}] ,\E)$ and $\M\times 
\dis\Pi_{i=1}^{n}\mathcal{R}^{0}([\t_{i-1},\t_{i}] ,\M)\times\mathcal{R}^0([\t_{i-1},\t_{i}] ,\E)$. This proof will be will be organized in six steps.
\end{com}

\bigskip

\so{\bf step 1}: {\it ${\mathcal{P}}(\t,\mathfrak{U})$ and (resp. ${\mathsf{P}}(\t,\mathfrak{U})$) is open for the compact-open topology on ${\mathcal{P}}(E)$ (resp. ${\mathcal{P}}(M)$)} (cf. $\S$ \ref{topologyPE}).\\
Since $\pi$ is open and $\pi({\mathcal{P}}(\t,\mathfrak{U}))={\mathsf{P}}(\t,\mathfrak{U})$, we have only to prove this result for ${\mathcal{P}}(\t,\mathfrak{U})$.
For any open set $O$ in $M$ and any subinterval $[\a,\b]$ of $I$, since $[\a,\b]$ is compact and $\pi\circ\gamma$ is continuous the set 
$$\{c\in \mathcal{P}(E)\;:\; \pi\circ c([\a,\b])\subset O\}$$
is an open subset of $\mathcal{P}(E)$. Thus the sets $\mathcal{O}_i= \{ c\in \mathcal{P}(E) \;|\; {\pi\circ c([\t_{i-1},\t_i])}\subset U_i\}$ are open and $ {\mathcal{P}}(\t,\mathfrak{U})$ is exactly the intersection $\bigcap_{i=1}^n\mathcal{O}_i$ which ends the proof.
 

\bigskip

\so{\bf step 2}: {\it $\mathcal{F}_{\t,\mathfrak{U}}$ and $\mathsf{f}_{\t,\mathfrak{U}}$ are injective and continuous.}\\
We will prove it only for $\mathcal{F}_{\t,\mathfrak{U}}$ as the proof for $\mathsf{f}_{\t,\mathfrak{U}}$ is analogous. Let $c$ and $c'$ be two paths in ${\mathcal{P}}(\mathfrak{U}, \t)$ such that $\mathcal{F}_{\t,\mathfrak{U}}(c)=\mathcal{F}_{\t,\mathfrak{U}}(c')$. According to the previous notations, this implies for $i=1$ that
$$\phi_1(c(0))=\phi_1(c'(0)),\; (\phi_1(\pi\circ c_1))^{(1)}=(\phi_1(\pi\circ c'_1))^{(1)},\; u_1=u'_1$$
Therefore $c_{| [\t_0,\t_1]}=c'_{| [\t_0,\t_1]}$.
Assume that for $1\leq i<n$ we have $c_{|[\t_0,\t_{i-1}]}=c'_{|[\t_0,\t_{i-1}]}$.
Then we can write $c_{| [\t_{i-1},\t_{i}]}=\Phi_i^{-1} (\hat{\g}_{i}, u_{i})$ and $c'_{| [\t_{i-1},\t_{i}]}=\Phi_{i}^{-1} (\hat{\g}'_{i}, u'_{i})$. Now from our assumption we have 
$\hat{\g}_{i}^{(1)}=(\hat{\g}')_{i}^{(1)}$ and $u_{i}=u'_{i}$. Since $c(\t_i)=c'(\t_i)$ so $\phi_i(\pi\circ c(\t_{i-1}))=\phi_i(\pi\circ c'(\t_{i-1}))$. But according to Remark \ref{restrictstrong} the restriction of $\hat{\g}_{i}^{(1)}$ (resp $(\hat{\g}')_{i}^{(1)}$) coincides with $\phi_i(\pi\circ c_{|[\t_{i-1},\t_i[})$ (resp. $\phi_i(\pi\circ c'_{|[\t_{i-1},\t_i[})$). Therefore this implies that $\hat{\g}_i=\hat{\g}'_i$ and finally $c=c'$ on $[\t_0,\t_i]$.
By induction on $1\leq i\leq n$ we obtain $c=c'$ and so $\mathcal{F}_{\t,\mathfrak{U}}$ is injective.

\bigskip

For the proof of continuity of $\mathcal{F}_{\t,\mathfrak{U}}$, it is sufficient to prove the continuity of each component of this map. 
The first component is ${\bf Ev^0_M}$ which is continuous as we have already seen (Proposition \ref{basecont} Point (3)). 
Now for $J=[\t_{i-1},\t_j]$ by application of Proposition \ref{basecont} Point (2), Proposition \ref{prop:restriction} and Proposition \ref{composition} we prove the $\tilde{\Phi}_i$ is continuous for $i=1,\dots,n$.\\

\so{\bf step 3}: {\it The range of $\mathcal{F}_{\t,\mathfrak{U}}$ (resp. $\mathsf{f}_{\t,\mathfrak{U}}$) is open in the Banach space\\ 
\indent $\M\times \dis\Pi_{i=1}^{n}\mathcal{R}^{0}([t_{i-1},t_{i}] ,\M)\times\mathcal{R}([t_{i-1},t_{i}] ,\E)$ (resp. $\M\times \dis\Pi_{i=1}^{n}\mathcal{R}^{0}([\t_{i-1},\t_{i}] ,\M)$).}\\
According to (\ref{n=1proj}), we have only to prove the result for $\mathcal{F}_{\t,\mathfrak{U}}$.
Recall that norm on each space $ \mathcal{R}^{0}([\t_{i-1},\t_{i}] ,\M)$ is
$$\nnorm{(x_i)}_\infty^0=\ \sup\{\norm{x_i(t)}_\M\; t\in[\t_{i-1},\t_{i}] \},$$
as defined in \ref{nota}.

{\bf Fix some $c\in {\mathcal{P}}(\t,\mathfrak{U})$.} We set $\bar{\g}_i=\phi_i\circ \pi \circ c_i$ and so we have 
$$\mathcal{F}_{\mathfrak{\t,U}}(c)=(\bar{\g}_1(\t_0), \bar{\g}_1^{(1)}, u_1\dots,\bar{\g}_{n}^{(1)}, u_{n}).$$
Since for $i=1,\dots, n-1$, each $\phi_i (U_i)$ and $\phi_i (U_i\cap U_{i+1})$ are open in $\M$, $\;[\t_{i-1},\t_{i}]$ is compact we can find $\d_i>0$ such that the open balls satisfy the following inclusions:

\begin{enumerate}
\item[(1)] $B(\bar{\g}_i(\t_{i}),\d_i)\subset \phi_i (U_i\cap U_{i+1})$ for $i=1,\dots,n-1$

\item[(2)] $B(\bar{\g}_i(t),\d_i)\subset \phi_i (U_i)$ for all $t\in [\t_{i-1},\t_{i}]$ and $i=1,\dots,n$.
\end{enumerate}


\noindent Moreover since $\phi_{i+1}\circ \phi_{i}^{-1}$ is continuous we can also find strictly positive numbers $\varepsilon_1,\dots,\varepsilon_{n-1}$ and $\eta_1,\dots,\eta_{n-1}$ such that:

\begin{enumerate}
\item[(3)] $\varepsilon_{i}<\d_i/2$, $\eta_{i}<\d_i/2$ for $i=1,\dots, n-1$
 
 
\item[(4)] $\phi_{i+1}\circ\phi_i^{-1}(B(\bar{\g}_i(\t_{i}), \d_i)\subset B(\bar{\g}_{i+1}(\t_i),\varepsilon_{i+1})$.
\end{enumerate}

\noindent We set $\varepsilon=\max(\varepsilon_1,\ldots \varepsilon_{n-1})$ and $\eta=\max(\eta_1,\ldots \eta_{n-1})$.

We denote by $\mathcal{O}$ the open neighbourhood of $\mathcal{F}_{\t,\mathfrak{U}}(c)$ consisting of $ (x',x'_1,u'_1,\dots,x'_{n},u'_{n})\in \M\times \dis\Pi_{i=1}^{n}\mathcal{R}^{0}([\t_{i-1},\t_{i}] ,\M)\times\mathcal{R}^0([\t_{i-1},\t_{i}] ,\E)$ such that 
$$ \norm{\bar{\g}(t_0)-x'}_\M<\varepsilon,\quad \nnorm{\bar{\g}^{(1)}_i-x'_i}_\infty^i<\eta,\qquad\qquad 1\leq i\leq n.$$

\noindent Now we will show that $\mathcal{O}$ is contained in $\mathcal{F}_{\t,\mathfrak{U}}(\mathcal{P}^k(\t,\mathfrak{U}))$.
To this end to any $(x',x'_1,u'_1,\dots,x'_n,u'_n)\in \mathcal{O}$ we associate a $1$-regulated path
$(\bar{\g}'_1,\dots,\bar{\g}'_n)$ in $ \dis\Pi_{i=1}^{n}\mathcal{R}^{1}([\t_{i-1},\t_{i}] ,\M)$ such that

\begin{enumerate}
\item[(5)] $\bar{\g}'_i([\t_{i-1},\t_{i}]\subset \phi_i(U_i),\; \bar{\g}'_{i+1}(\t_i)=\phi_{i+1}\circ\phi_i^{-1}(\bar{\g}'_{i}(\t_i)) \textrm{ and } \bar{\g}'_i(\t_i)\in B(\bar{\g}_{i}(\t_i),\d_i)$
\end{enumerate}

\noindent in the following way:

$\bullet$ for $i=1$ we define $\bar{\g}'_1$ by $\bar{\g}'_1=x'+\dis\int_{\t_0}^t x'_1(s)ds$. By construction $\bar{\g}'_1$ is $1$-regulated and from the property (3) we have 
$\norm{\bar{\g}_1-\bar{\g}'_1}_\infty<\d_1$. This implies that in particular $\bar{\g}'_1(\t_1)$ belongs to $B(\bar{\g}_1(\t_1),\d_1)$ and from the property (2) we get 
$\bar{\g}'_1([\t_{0},\t_{1}])\subset \phi_1(U_1)$;\\

$\bullet$ assume that we have built a sequence of $1$-regulated paths 
$(\bar{\g}'_1,\dots,\bar{\g}'_i)$ in $ \dis\Pi_{j=1}^{i}\mathcal{R}^{0}([\t_{j-1},\t_{j}] ,\M)$ which satisfies the property (5). 
Since $\bar{\g}'_i(\t_i)\in B(\bar{\g}_i(\t_i),\d_i)$ then from property (4), $\phi_{i+1}\circ\phi_i^{-1}(\bar{\g}'_i(\t_i))$ belongs to $B(\bar{\g}_{i+1}(\t_i),\varepsilon_{i+1})$.
We set 
$$\bar{\g}'_{i+1}=\phi_{i+1}\circ\phi_i^{-1}(\bar{\g}'_i(\t_i))+\dis\int_{\t_i}^t x'_{i+1}(s)ds.$$

\noindent Therefore we get $\bar{\g}'_{i+1}(\t_i)=\phi_{i+1}\circ\phi_i^{-1}(\bar{\g}'_i(\t_i))$ and again as for $i=1$ we also have 
$\bar{\g}'_{i+1}([\t_{i},\t_{i+1}])\subset \phi_{i+1}(U_{i+1})$.
We can now define a path $c'\in \mathcal{P}(\t,\mathfrak{U})$ such that 
$\mathcal{F}_{\t,\mathfrak{U}}(c')=(x',x'_1,u'_1,\dots,x'_n,u'_n)$ by setting
$$c'_{| [\t_{i-1},\t_i]}=\Phi_i^{-1}(\bar{\g}'_i, u'_i).$$
Thus we have proved that $\mathcal{O}$ is contained in $\mathcal{F}_{\t,\mathfrak{U}}(\mathcal{P}^k(\t,\mathfrak{U}))$.

\smallskip

\so{\bf step 4}:{ \it Consider two partitions $\t=(\t_i)_{i=0,\dots,n}$ and $\s=(\s_j)_{j=0,\dots,m}$ of $I$ with associated sequences $\mathfrak{U}=(U_i)_{i=1,\dots,n}$ and $\mathfrak{V}=(V_j)_{j=1,\dots,m}$ of domains of charts $(U_i,\phi_i)$ and $(V_j,\psi_j)$ respectively such that $\mathcal{P}(\t,\mathfrak{U})\cap \mathcal{P}(\s,\mathfrak{V})\not= \emptyset$. Then 
the map 
 $\mathcal{F}_{\t\mathfrak{U}}\circ (\mathcal{F}_{\s,\mathfrak{V}})^{-1}$ (resp. $(\mathsf{f}_{\t,\mathfrak{U}}\circ (\mathsf{f}_{\s\mathfrak{V}})^{-1}$) is a diffeomorphism from it definition domain onto its range.} \\

\noindent 
We will use the ``reconstruction map''

\noindent $\Theta: \M\times \dis\Pi_{i=1}^n\mathcal{R}^0([\t_{i-1},\t_i],\M)\ap \dis\Pi_{i=1}^n\mathcal{R}^{1}([\t_{i-1},\t_i],\M)$ defined by

$$\Theta(x,y_1,\dots,y_n)=(\bar{\g}_1,\dots\bar{\g}_n),$$
where 
\begin{align*}
\bar{\g}_1 &= x+\dis\int_{\t_0}^t y_{1}(s)ds,\\
\bar{\g}_{i+1} &= \phi_{i+1}\circ\phi_i^{-1}(\bar{\g}_i(\t_i))+\dis\int_{\t_i}^t y_{i+1}(s)ds\qquad\textrm{ for }i=1,\dots, n-1.
\end{align*}
We need the following property:
 

\begin{lem}\label{smoothpsi} $\Theta$ is a smooth map.
\end{lem}
\proof{\it Proof of Lemma \ref{smoothpsi}}\\
It is sufficient to prove that each component $\Theta_i:(x,y_1,\dots,y_n)\mapsto \bar{\g}_i$ is smooth. For any Banach space $\mathbb{F}$, denote again by $ {\bf Ev^t_\mathbb{F}}$ the evaluation map form any Banach space $\mathcal{R}^k(I,\mathbb{F})$ to $\F$, for $k=0, 1$, given by ${ \bf Ev^t_\mathbb{F}}(c)=c(t)$. Since this map is continuous linear, so it is smooth. Now for $2\leq i\leq n$ we can write 
$$\Theta_i={\bf Ev^{\bf \t_{i-1}}}(\phi_{i+1}\circ\phi_i^{-1}(\bar{\g}_i))+L(y_i)$$
where $L$ is linear map from $\mathcal{R}^0([\t_{i-1},\t_i],\M)$ to $\mathcal{R}^{1}([\t_{i-1},\t_i],\M)$ and since $\t_i-\t_{i-1}<1$ then $L$ is a contraction and so $L$ is continuous which implies that $L$ is smooth. 
 The same type of arguments can be applied to $\Theta_1$.\\ 
 \endproof\\
 \bigskip
 \noindent {\bf We come back to the proof of this step}. We consider the map
 $\mathcal{F}_{\t,\mathfrak{U}}\circ (\mathcal{F}_{\s,\mathfrak{V}})^{-1}$ from $\mathcal{F}_{\s,\mathfrak{V}}(\mathcal{P}(\t,\mathfrak{U})\cap \mathcal{P}(\s, \mathfrak{V}))\subset \M\times \dis\Pi_{j=1}^m\mathcal{R}^0([\s_{j-1},\s_j],\M)\times\mathcal{R}([\s_{j-1},\s_j],\E)$ to $\M\times \dis\Pi_{i=1}^n\mathcal{R}([\t_{i-1},\t_i],\M)\times\mathcal{R}([\t_{i-1},\t_i],\E)$.
In fact $\mathcal{F}_{\s,\mathfrak{V}}(\mathcal{P}(\t,\mathfrak{U})\cap \mathcal{P}(\s,\mathfrak{V}))=\mathsf{f}_{\s,\mathfrak{V}}(\mathsf{P}(\t,\mathfrak{U})\cap \mathsf{P}(\s,\mathfrak{V}))\times\Pi_{i=1}^n \mathcal{R}([\s_{j-1},\s_j],\E)$ and the range of $\mathcal{F}_{\t,\mathfrak{U}}\circ (\mathcal{F}_{\s,\mathfrak{V}})^{-1}$ is 
$\mathsf{f}_{\t,\mathfrak{U}}(\mathsf{P}(\t,\mathfrak{U})\cap \mathsf{P}(\s,\mathfrak{V}))\times\Pi_{i=1}^n \mathcal{R}([\t_{j-1},\t_j],\E)$.

For $c\in {\mathcal{P}}(\s,\mathfrak{V})$ we set $\hat{\g}_i=\psi_i\circ \pi \circ c_i$ and so we have 
\begin{equation}\label{FsV}
\mathcal{F}_{\s,\mathfrak{V}}(c)=(\hat{\g}_1(0), \hat{\g}_1^{(1)}, v_1,\dots,\hat{\g}_{m}^{(1)}, v_{m}).
\end{equation}

Let $\hat{\Theta}$ be the reconstruction map associated to $(\mathfrak{V},\s)$ and in order to make further formulas more readable we will use notation $\hat{\Theta}(\mathfrak{w})$ in the following sense: $\mathfrak{w}=(x,y_1,v_1,\dots, y_m, v_m)\mapsto \hat{\Theta}(\mathfrak{w})=(\hat{\g}_1,\dots,\hat{\g}_m)$.

 \noindent Now the map:
 $$ \mathcal{F}_{\t,\mathfrak{U}}\circ (\mathcal{F}_{\s,\mathfrak{V}})^{-1}: \mathsf{f}_{\s,\mathfrak{V}}(\mathsf{P}(\t,\mathfrak{U})\cap \mathsf{P}(\s,\mathfrak{V}))\times\Pi_{i=1}^n \mathcal{R}([\s_{j-1},\s_j],\E)\ap \M\times \dis\Pi_{i=1}^n\mathcal{R}^0([\t_{i-1},\t_i],\M)\times\mathcal{R}([\t_{i-1},\t_i],\E) $$
can be written as:
$$ \mathcal{F}_{\t,\mathfrak{U}}\circ (\mathcal{F}_{\s,\mathfrak{V}})^{-1}(\mathfrak{w})=(\mathcal{G}_0(\mathfrak{w}), \mathcal{G}_1(\mathfrak{w}),\mathcal{W}_1(\mathfrak{w})\cdots,\mathcal{G}_n(\mathfrak{w}),\mathcal{W}_n(\mathfrak{w})).$$
 
To prove that $\mathcal{F}_{\t,\mathfrak{U}}\circ(\mathcal{F}_{\s,\mathfrak{V}})^{-1}$ is smooth it is sufficient to prove all the components $\mathcal{G}_0$ and $(\mathcal{G}_i, \mathcal{W}_i)$ for $i=1,\dots,n$ are smooth. The first one is simply
a map $\mathfrak w \mapsto x\mapsto \phi_1\circ \psi_1^{-1} (x)$. It is smooth since it is composition of projection onto first term and transition map of the manifold $M$.
The proof for the other components is not as so easy.

For $i=1,\dots,n$ and $j=1,\dots,m$ we denote by $I_{ij}$ the intersection $[\t_{i-1},\t_i]\cap [\s_{j-1},\s_j]$ and we denote by $N_{\t_i}=\{j: I_{ij}\not=\emptyset\}$ and 
$N_{\s_j}=\{i: I_{ij}\not=\emptyset\}$. 
Thus we have $[\s_{j-1},\s_{j}]=\dis\cup_{i\in N_{\s_j}} I_{ij}$ and $[\t_{i-1},\t_{i}]=\dis\cup_{j\in N_{\t_i}} I_{ij}$. From notations in relation (\ref{FsV}), we consider the restrictions\footnote{in the sense of Definition \ref{restrict}}
$\hat{\g}_{ij}$ 
of $\hat{\g}_j$ to $I_{ij}$ and $\hat{\g}_j(I_{ij})$ is contained in $\psi_j(U_i\cap V_j)$. In the same way, if $v_{ij}$ is the restriction of $v_j$ to $I_{ij}$ (in the usual sense), then $w_j(I_{ij})$ is contained in $\psi_j(U_i\cap V_j)\times \E$.\\

Now for $i=1,\dots,n$, the pair of components $(\mathcal{G}_i, \mathcal{W}_i)$ is the composition of the following smooth maps :
\begin{description}
 \item[(I)] ${r}_{ij}: \M\times \dis\Pi_{j=1}^m\mathcal{R}^0([\s_{j-1},\s_j],\M)\times\mathcal{R}^0([\s_{j-1},\s_j],\E) \ap \mathcal{R}^{1}(I_{ij}, \M)\times \mathcal{R}^0(I_{ij}, \E)$ defined by 

$r_{ij}(\mathfrak{w})=(\mathsf{r}_{I_{ij}}\circ\hat\Theta_j(\mathfrak w), v_{ij})$. This map is smooth as the composition of the continuous linear map $(\mathsf{r}_{I_{ij}},\textrm{ restriction to } I_{ij})$
with the map $(\hat{\Theta}_j, pr_j)$ which is smooth by Lemma \ref{smoothpsi} and where $pr_j$ is the projection of $\M\times \dis\Pi_{j=1}^m\mathcal{R}^0([\s_{j-1},\s_j],\M)\times\mathcal{R}^0([\s_{j-1},\s_j],\E)$ onto $ \mathcal{R}^0([\s_{j-1},\s_j], \E)$.\\ 

\item[(II)] $d_{ij}: \mathcal{R}^{1}(I_{ij}, \M)\times \mathcal{R}^0(I_{ij}, \E) \ap \mathcal{R}^{0}(I_{ij}, \M)\times \mathcal{R}^0(I_{ij}, \E)$ defined by

$d_{ij}(f, v_{ij})=\big((\phi_i\circ \psi_j^{-1}\circ f)^{(1)}, g_{ij}(\phi_i\circ \psi_j^{-1}\circ \g)v_{ij}\big)$ where $g_{ij}$ is the smooth field of automorphisms of $\E$ associated to $\Phi_i\circ \Psi_j^{-1}\equiv(\phi_i\circ \psi_j^{-1}, g_{ij})$. The smoothness of the map $d_{ij}$ is a consequence of the Proposition \ref{composition} and the smoothness of $g_{ij}$ in the two variables.

\item[(III)] $C_i:\dis\Pi_{j\in N_{\t_i}} \mathcal{R}^{0}(I_{ij}, \M)\times \mathcal{R}^0(I_{ij}, \E)\ap \mathcal{R}^0([\t_{i-1},\t_i], \M)\times\mathcal{R}^0([\t_{i-1},\t_i], \E)$ defined by

$C_i((u_{ij},v_{ij})_{j\in N_{\t_i}})=( u_i,v_i)$ where $u_i$ is the concatenation of the family $(u_{ij})_{j\in N_{\t_i}}$ and $(v_i)_{| I_{ij}}=v_{ij}$. This map is clearly linear.

If we provide $\mathcal{R}^{0}(I_{ij}, \M)\times\mathcal{R}^0(I_{ij},\M)$ with the norm
$$\nnorm{(u_{ij},v_{ij})}_\infty^{ij}= \sup_{t\in I_{ij}} \norm{(u_{ij}(t), v_{ij}(t))}_{\M\times\E}$$

then $\dis\sum_{j\in N_{\t_i}} \nnorm{\;}_\infty^{ij}$ is a norm on 
$\mathcal{R}^0([\t_{i-1},\t_i], \M)\times\mathcal{R}^0([\t_{i-1},\t_i], \E)$ and on\\ $\dis\Pi_{j\in N_{\t_i}} \mathcal{R}^{0}(I_{ij}, \M)\times\mathcal{R}^{0}(I_{ij}, \E)$. 
In this 
way $C_i$ is an isometry and so it is smooth.\\
\end{description}

\noindent The result concerning $(\mathsf{P}(\mathfrak{U}, \t),\mathsf{f}_{\t,\mathfrak{U}})$ is implied by application of the previous arguments to the subspace $ \mathsf{f}_{\mathfrak{V},\s}(\mathsf{P}(\t,\mathfrak{U})\cap \mathsf{P}(\s,\mathfrak{V}))\times\{0\}$ of 
$ \mathsf{f}_{\mathfrak{V},\s}(\mathsf{P}(\t,\mathfrak{U})\cap \mathsf{P}(\s,\mathfrak{V}))\times\Pi_{i=1}^n \mathcal{R}^0([\s_{j-1},\s_j],\E)$.\\

\noindent From all the properties established in the whole of the previous steps, it follows that the set of all the possible pairs $(\mathcal{P}(\t,\mathfrak{U}),\mathcal{I}_\t\circ\mathcal{F}_{\s,\mathfrak{V}})$ (resp. $(\mathsf{P}(\t,\mathfrak{U}),\mathsf{I}_\t\circ \mathsf{f}_{\s,\mathfrak{V}})$) gives rise to an atlas on ${\mathcal{P}}(E)$ (resp. ${\mathcal{P}}(M)$) whose any chart domain is open for the topology induced by ${\mathcal{P}}(E)$ (resp. ${\mathcal{P}}(M)$). Since $\mathcal{P}(E)$ (resp. $\mathsf{P}(M)$) is a Hausdorff topological space it follows that we provide a Banach structure.

\bigskip

\so{\bf step 5} : {\it $\pi: \mathcal{P}(E)\ap \mathsf{P}(M)$ is a Banach bundle.}

Now on one hand from the initial definition at the beginning of the proof we have the following diagram
 $$\begin{tikzcd}
\mathcal{P}(\t,\mathfrak{U})\arrow{r}{\mathcal{F}_{\t,\mathfrak{U}}}\arrow{d}{}
&\mathsf{f}_{\t,\mathfrak{U}}(\mathsf{P}(\t,\mathfrak{U}))\times \dis\Pi_{i=1}^n\mathcal{R}^0([\t_{i-1},\t_i],\E) \arrow{d}{q_1}\\
\mathsf{P}(\t,\mathfrak{U})\arrow{r}{\mathsf{f}_{\t,\mathfrak{U}}}&\mathsf{f}_{\t,\mathfrak{U}}(\mathsf{P}(\t,\mathfrak{U}))
\end{tikzcd}
$$

On the other hand from step 4 recall that we have that
the map $(\mathcal{F}_{\t,\mathfrak{U}})\circ (\mathcal{F}_{\s,\mathfrak{V}})^{-1}$ is a map from $\mathsf{f}_{\mathfrak{V},\s}(\mathsf{P}(\t,\mathfrak{U})\cap \mathsf{P}(\s,
\mathfrak{V})))\times\Pi_{j=1}^n \mathcal{R}^0([\s_{j-1},\s_j],\E)$ to 
$\mathsf{f}_{\t,\mathfrak{U}}(\mathsf{P}(\t,\mathfrak{U})\cap \mathsf{P}(\s,\mathfrak{V}))\times\Pi_{i=1}^n \mathcal{R}^0([\t_{j-1},\t_j],\E)$. \\ 

From the previous diagram, and according to the proof of step 4 we have:

\begin{description} 
\item[-] if ${\g}\in \mathsf{f}_{\mathfrak{V},\s}(\mathsf{P}(\t,\mathfrak{U})\cap \mathsf{P}(\s,\mathfrak{V}))$ this map sends each factor of type $\{\mathsf{f}_{\mathfrak{V},\s}(\mathsf{P}(\t,\mathfrak{U})(\g)\}\times\mathcal{R}^0(I_{ij}, \E)$ onto the factor $\{\mathsf{f}_{\t,\mathfrak{U}}(\mathsf{P}(\t,\mathfrak{U})(\g)\}\times\mathcal{R}^0(I_{ij}, \E)$

 \item[-] we apply point {\bf (II)} in restriction to this factor
 
 \item[-] we compose with $C_i$ (cf. Point {\bf (III)}) in restriction to the same factor. 
 \end{description} 
 
 \noindent In this way, we see the map $(\mathcal{F}_{\t,\mathfrak{U}})\circ (\mathcal{F}_{\s,\mathfrak{V}})^{-1}$ is a map from \\
$\mathsf{f}_{\mathfrak{V},\s}(\mathsf{P}(\t,\mathfrak{U})\cap \mathsf{P}(\s,\mathfrak{V}))\times\Pi_{j=1}^m \mathcal{R}^0([\s_{j-1},\s_j],\E)$ to 
$\mathsf{f}_{\t,\mathfrak{U}}(\mathsf{P}(\t,\mathfrak{U})\cap \mathsf{P}(\s,\mathfrak{V}))\times\Pi_{i=1}^n \mathcal{R}^0([\t_{j-1},\t_j],\E)$. of type 
$({\bf z,w)}\mapsto ( \mathsf{f}_{\t,\mathfrak{U}}\circ\mathsf{f}_{\mathfrak{V}\s}^{-1} (\bf{z}), G( \mathsf{f}_{\t,\mathfrak{U}}\circ\mathsf{f}_{\mathfrak{V}\s}^{-1} (z)){\bf w})$ where ${\bf G(z')}$ belongs to\\ $GL(\dis\Pi_{i=1}^n\mathcal{R}^0([\s_{i-1},\s_i],\E),\dis\Pi_{i=1}^n\mathcal{R}^0([\t_{i-1},\t_i],\E)$ and is smooth relative to ${\bf z}'$.

\bigskip

\so{\bf step 6} : {\it existence of Banach isomorphisms $\mathcal{T}_c: T_c\mathcal{P}(E)\ap \G(c)$ and $\mathsf{T}_\g: T_\g\mathsf{P}(M)\ap \G(\g)$ }

\noindent Similarly as in previous steps, the result for $\mathsf{T}_\g: T_\g\mathsf{P}(M)\ap \G(\g)$ is implied by result for $\mathcal{T}_c: T_c\mathcal{P}(E)\ap \G(c)$ since $\mathsf{T}_\g=T\pi\circ \mathcal{T}_c$ if $\g=\pi\circ c$.\\
Now any $V\in T_c {\mathcal{P}}(E)$ is defined by a $C^1$ path $\tilde{c}:[-\d,\d]\ap {\mathcal{P}}(E)$ (called a {\it $C^1$-deformation of $c$}) such that 
$$\dis\frac{d \tilde{c}}{d\epsilon}_{| \epsilon=0}=V$$
To this path $\tilde{c}$ is canonically associated a map from $[-\d,\d]\times [0,1]$ to $\mathcal{P}(E)$ also denoted $\tilde{c}$.
 Consider a chart $(\mathcal{P}(\t,\mathfrak{U}), \mathcal{F}_{\t, \mathfrak{U}})$ around $c$. 
For $\d$ small enough, for any $i=0,\dots,n-1$, the restriction $\tilde{c}_i$ of $\tilde{c}$ to $[\t_{i-1},\t_{i}]$ must be contained in $E_{| U_i}$. The evaluation from $\mathcal{P}_{[\t_{i-1},\t_{i}]}(\phi_i(U_i)\times E)$ to $\M\times\E$ is denoted ${\bf Ev^t_{E_{| U_i}}}$. We have
${\bf Ev^t_{E_{| U_i}}}\circ \mathcal{F}_{\t,\mathfrak{U}}(c)=\Phi_i(c_i)(t)$ if $t\in [\t_i,\t_{i+1}]$.

\noindent The evaluation map is linear and smooth and it commutes with differential. Therefore we have
${\bf Ev^t_{E_{| U_i}}}\circ \mathcal{F}_{\t,\mathfrak{U}}(\tilde{c}(\epsilon))=\Phi_i(\tilde{c}_i(\epsilon))(t)=\Phi_i\circ {\bf Ev^t}(\tilde{c}(\epsilon))$ if $t\in [\t_i,\t_{i+1}]$ and 
$\{{\bf Ev^t_{E_{| U_i}}}\circ \mathcal{F}_{\t,\mathfrak{U}}(\tilde{c})\}_{|\epsilon=0}=\dis\frac{\p \Phi_i\circ \tilde{c}}{\p \epsilon}(0,t)=D\Phi_i\circ \dis\frac{\p \tilde{c}}{\p \epsilon}(0,t)={\bf Ev^t_{E_{| U_i}}}\circ D\Phi_i(V_i(t))$ where $V_i=V_{| [\t_{i-1},\t_i]}$. Now the usual transformation formulae between charts in the tangent bundle of a Banach manifold imply that $V$ is a well defined section of $\G(c)$. \\
Conversely, if $V\in T_c\mathcal{P}(E)$, with the previous notations each $V_i$ belongs to $\G(c_i)$ and so can be written $(\bar{\varphi}_i,\bar{\vartheta}_i )$ via the representation of $\mathcal{T}^{\Phi_i}_{c_i}$ of $\G(c_i)$ with $\mathcal{R}^1([\t_{i-1},\t_i] ,\M)\times \mathcal{R}^0([\t_{i-1},\t_i] ,\E)$ induced by $\Phi_i$. Now fix some representation $\mathcal{T}_\g$ and $\mathcal{T}_c$ of $\G(\g)$ and $\G(c)$ ({cf.} Notations and conventions \ref{triviaVF})). According to Remark \ref{localtrivial}, we have a $1$-regulated (resp. $0$-regulated) field $A_{\g_i}$ (resp. $A_{c_i}$) of automorphisms of $\M$ (resp. $\M\times\E$ ) such that 
$$\mathcal{T}_{\g_i}^{\phi_i}=A_{\g_i}\circ \mathcal{T}_{\g_i} ,\;\; \mathcal{T}_{c_i}^{\Phi_i}=A_{c_i}\circ \mathcal{T}_{c_i} \textrm{ and } A_{\g_i}=q_1\circ A_{c_i}.$$
where $\mathcal{T}^{\phi_i}_{\g_i}$ is the representation of $\G(\g_i)$ associated to $\phi_i$. Therefore if $(\varphi_i,\vartheta_i)=A_{c_i}^{-1}((\bar{\varphi}_i,\bar{\vartheta}_i )$ we also have $\varphi_i=A_{\g_i}^{-1} (\bar{\varphi}_i)$ and moreover the pairs $\{(\varphi_i,\vartheta_i)\}_{i=1,\dots,n}$ stick together into a global $(\varphi,\vartheta)$ such that the restriction of $(\varphi,\vartheta)$ to $[\t_{i-1},\t_i]$ is precisely $(\varphi_i,\vartheta_i)$. Clearly this implies that $\mathcal{T}_{c_i}^{-1}(\varphi,\vartheta)=V$ and moreover $\mathcal{T}_\g^{-1}(\varphi)=T\pi(V)$ which ends the proof.\\
\endproof

\appendix

\section{Banach manifold structure in the case of $C^1$ curves}
{\it In this Appendix we give an adapted version of the results of previous sections for curves $c\in\mathcal{P}(E)$ which are $C^1$ and whose projection on $M$ is $C^2$. This context is explicitly used \cite{Pe}.}

\begin{rem-a} \label{ck-1gk1}\normalfont Let us consider the spaces $\mathsf{C}^{k-1}(M)$ of curves $\g:I\ap M$ of class $C^k$ and $\mathcal{C}^k(E)$ of lifts to $E$ of curves in $\mathsf{C}^{k-1}(M)$ which are of class $C^{k-1}$. Then for any $\g\in \mathsf{C}^{k-1}(M)$ by same arguments as in subsection \ref{GcGg} (via trivializations $\mathcal{T}_\g$ and $\mathcal{S}_\g$) we can define on the set $\G^k(\g)$ of vector fields of class $C^k$ along $\g$ a structure of Banach space modeled on $\mathsf{C}^k(\M)$. We can also define on the set $\G^{k-1}(c)$ of vector fields of class $C^{k-1}$ along a lift $c\in \mathcal{C}^k(E)$ a structure of Banach space modeled on $\mathsf{C}^k(\M)\times \mathsf{C}^{k-1}(\E)$.

More precisely for $k=1$ all the arguments used in the subsection \ref{GcGg} are valid if we replace ``$0$-regulated'' by ``$C^0$'' and ``$1$-regulated'' by ``$C^1$''. For $k>1$, on one hand all the arguments before Proposition \ref{Banachstr} are valid by replacing ``regulated'' by ``$C^{k-1}$'' and ``$1$-regulated'' by ``$C^k$''. On the other hand, the proof of a corresponding version of Proposition \ref{Banachstr} is obtained by using the expression of the transformation of the $l$-order derivative of a curve under transformation $\mathcal{T}_\g$ or $\mathcal{S}_\g$ for $1\leq l\leq k$ and comparable arguments concerning bounds (as for $k=1$) relative to all the semi-norm $\sup_{t\in I} \norm{\d^{(l)}(t)}_\F$ for any curve $\d: I\ap \F$ of class $C^l$ in a Banach space $\mathbb{F}$.
\end{rem-a}

\begin{rem-a}\label{ck-1gk2}
\normalfont 
For any $l\geq 1$ and any Banach manifold $B$ let $T_x^l B$ be the set of equivalence classes of local germs of curves which have a contact at order $l$ at $x$. Then $T^lB=\dis\cup_{x\in N}T^l_xB$ has a structure of Banach manifold modeled on $\mathbb{B}^{l+1}$ (see \cite{Su} Theorem 2.1). Moreover, $p_B^l:T^l B\ap B$ is the natural projection $(p_B^l)^{-1}$ is a chart domain for each domain chart $U$ in $B$\footnote{ $p_B^l:T^lB\ap B$ is a Banach bundle if and only if there exists a linear connection on $TB$ (see \cite{Su} Theorem 2.11)}. Note that $T^1B=TB$ and for simplicity we set $T^0B=B$.
 
As previously, we can provide $\mathcal{C}^k(E)$ with the topology generated by sets $\mathcal{N}(K,V_1,\cdots,V_{k},W_1,\cdots, W_{k+1})$ of paths $ c \in \mathcal{P}(E) $ such that the closure of $c^{l}(K)$ is contained in $ V_l$, and $(\pi \circ c)^{(l)}(K)$ is contained in $W_l$ where $K$ is a compact subset of $I$ and $V_l$ is an open set in $T^l E$ (resp. $W_l$ an open in $T^lM$) for $l=0,\dots k$ (resp. for $l=0,\dots, k+1$). In the same way, we can define a topology on $\mathsf{C}^k(M)$. Again these topologies are Hausdorff and we have a natural projection $\pi^k:\mathcal{C}^k(E)$ onto $\mathsf{C}^k(M)$.
\end{rem-a}

\begin{rem-a}\label{ck-1gk3} \normalfont The Proposition \ref{basecont} remains true if we replace $\mathcal{P}$ by $\mathcal{C}^k$, $\mathcal{R}^1(I,\phi(U))$ by $C^k(I,\phi(U))$, $\mathcal{R}^0(I,\E)$ by $C^{k-1}(I,\E)$.
\end{rem-a}

\subsection{Banach structure on $\mathcal{C}(E)$}

To simplify the notation we denote by $\mathcal{C}(E)$ the set $\mathcal{C}^1(E)$ 
(curves $c\in \mathcal{P}(E)$ of class $C^1$ whose projection $\g=\pi\circ c$ on $M$ is of class $C^2$) 
and by $\mathsf{C}(M)$ the set $\mathsf{C}^1(M)$ 
(curves $\g\in \mathsf{P}(M)$ of class $C^2$). We provide $\mathcal{C}(E)$ and $\mathsf{C}(M)$ with the topology defined in Remark \ref{ck-1gk2}.

\bigskip

At first according to notation in Remark \ref{ck-1gk1} we have a version of result from section \ref{PEn}:

\begin{theo-a} \label{Cbanachn}${}$
\begin{enumerate}
\item The set $\mathsf{C}(M)$ has a Banach manifold structure modeled on the Banach space $C^2(I,\M)$ and any chart domain of this structure is open for its compact open topology. Moreover, we have a Banach isomorphism from $T_\g\mathsf{C}(M)$ onto $\G^2(\g)$. 
\item The set ${\mathcal{C}}(E)$ has a Banach manifold structure modeled on the Banach space $C^k(I,\M)\times C^{1}(I ,\E)$ such that each chart domain is open for the natural topology of ${\mathcal{C}}(E)$ and we have a Banach isomorphism from $T_c {\mathcal{C}}(E)$ onto $\G^2(c)$. Moreover $\pi:\mathcal{P}(E)\ap \mathsf{P}(M)$ is a Banach bundle with typical fiber $\mathsf{C}^1(I ,\E)$.
\end{enumerate}
\end{theo-a}
 
\bigskip

Applying Remark \ref{ck-1gk1}, Remark \ref{ck-1gk2} and Remark \ref{ck-1gk3} for $k=1$, it is easy to see that Lemma \ref{isoTCGc} is also valid in this context by replacing the spaces $\mathcal{R}^0([\t_{i-1},\t_i],\M)$, $\mathcal{R}([\t_{i-1},\t_i],\E)$
and $ \mathcal{R}^{1}(I,\M)$ respectively by $C^1([\t_{i-1},\t_i],\M)$, $C^1([\t_{i-1},\t_i],\E)$ and $C^{2}(I,\M)$. Then the proof of Theorem \ref{Cbanachn} is an adaptation, point by point, of the proof of Theorem \ref{banachP}.

\begin{rem-a}\label{Kri}\normalfont
Since $C^1([\t_{i-1},\t_i],\M)$, $C^1([\t_{i-1},\t_i],\E)$ and $C^{2}(I,\M)$ satisfy all the axioms in section 2 of \cite{Kri}, Theorem \ref{Cbanachn} can be also proved by applying with the same procedure but using directly arguments of \cite{Kri}.
\end{rem-a}

\section*{Acknowledgement}
This research was partially supported by National Science Centre, Poland/Fonds zur Förderung der wissenschaftlichen Forschung grant
``Banach Poisson-Lie groups and integrable systems'' number 2020/01/Y/ST1/00123.


\begin{thebibliography}{GMV15}

\bibitem[AMR02]{AMR}
R.~Abraham, J.~E. Marsden, T.~S. Ratiu: Manifolds, Tensor Analysis, and
  Applications.
\newblock Springer-Verlag, Berlin-Heidelberg, third edition, 2002.

\bibitem[AOP01]{AOP}
M.~Alcheikh, P.~Orro, F.~Pelletier: Singularit{\'e} de l’application
  extr{\'e}mit{\'e} pour les chemins horizontaux en g{\'e}om{\'e}trie
  sous-riemannienne.
\newblock \emph{Groupe de travail: Singularit{\'e}s et g{\'e}om{\'e}trie sous
  riemannienne, Chambéry octobre 1997}, 2001.

\bibitem[Arg20]{Ar}
S.~Arguillere: Sub-riemannian geometry and geodesics in {B}anach manifolds.
\newblock \emph{Journal of Geometric Analysis}, \textbf{30}(3):2897--2938,
  2020.

\bibitem[AS04]{AgSa}
A.~A. Agrachev, Y.~L. Sachkov: Control theory from the geometric viewpoint,
  \emph{Encyclopaedia of Mathematical Sciences}, volume~87.
\newblock Springer, 2004.

\bibitem[Die60]{Di}
J.~Dieudonn{\'e}: Foundations of modern analysis.
\newblock Academic Press, 1960.

\bibitem[Gl{\"o}15]{Glo}
H.~Gl{\"o}ckner: Measurable regularity properties of infinite-dimensional {L}ie
  groups.
\newblock \emph{arXiv preprint arXiv:1601.02568}, 2015.

\bibitem[GMV15]{GMV}
E.~Grong, I.~Markina, A.~Vasil’ev: Sub-{R}iemannian geometry on
  infinite-dimensional manifolds.
\newblock \emph{Journal of Geometric Analysis}, \textbf{25}(4):2474--2515,
  2015.

\bibitem[KM97]{KrMi}
A.~Kriegl, P.~W. Michor: The convenient setting of global analysis, volume~53.
\newblock American Mathematical Society, 1997.

\bibitem[Kri72]{Kri}
N.~Krikorian: Differentiable structures on function spaces.
\newblock \emph{Transactions of the American Mathematical Society},
  \textbf{171}:67--82, 1972.

\bibitem[Mic80]{Mic}
P.~W. Michor: Manifolds of differentiable mappings, \emph{Shiva mathematics
  series}, volume~3.
\newblock Shiva Publishing, 1980.

\bibitem[Mon02]{Mon}
R.~Montgomery: A tour of subriemannian geometries, their geodesics and
  applications, \emph{Mathematical Surveys and Monographs}, volume~91.
\newblock American Mathematical Society, 2002.

\bibitem[Pel20]{Pe}
F.~Pelletier: Homotopy classes and foliation in the {B}anach manifold of
  admissible curves for a {B}anach pre-{L}ie algebroid.
\newblock \emph{preprint LAMA june 2020}.

\bibitem[Pen67]{Pen1}
J.-P. Penot: Vari\'et\'es diff\'erentiables d'applications et de chemins.
\newblock \emph{C. R. Acad. Sci. Paris S\'er. A-B},
  \textbf{264}(24):A1066--A1068, 1967.

\bibitem[Pen12]{Pen2}
J.-P. Penot: Calculus without derivatives, \emph{Graduate Text in Mathematics},
  volume 266.
\newblock Springer, 2012.

\bibitem[PT01]{PiTa}
P.~Piccione, D.~V. Tausk: Variational aspects of the geodesics problem in
  sub-{R}iemannian geometry.
\newblock \emph{Journal of Geometry and Physics}, \textbf{39}(3):183--206,
  2001.

\bibitem[Sch16]{Sch}
A.~Schmeding: Manifolds of absolutely continuous curves and the square root
  velocity framework.
\newblock \emph{arXiv preprint arXiv:1612.02604}, 2016.

\bibitem[Sus73]{Su}
H.~J. Sussmann: Orbits of families of vector fields and integrability of
  distributions.
\newblock \emph{Transactions of the American Mathematical Society},
  \textbf{180}:171--188, 1973.

\end{thebibliography}

\end{document}